\documentclass{elsarticle}

\usepackage{amsfonts}
\usepackage{textcomp}
           
\usepackage{tikz}
\usetikzlibrary{arrows}

\usepackage{amssymb}
\usepackage{eqnarray} 

\usepackage{amsmath}
\usepackage{mathalfa} 
\usepackage{hyperref} 
\usepackage{pdfpages}		
\usepackage{dsfont}   		
\usepackage{algorithm}
\usepackage{algpseudocode}

\usepackage{bbm}
\usepackage{comment}

\usepackage{graphicx}
\usepackage[small]{caption}
\usepackage[utf8]{inputenc}
\usepackage{psfrag}
\usepackage{amsthm}
\usepackage{amsmath}
\usepackage{amssymb}
\usepackage{mathtools}

\usepackage{psfrag}
\usepackage{url}

\newcommand{\sD}{\mathcal{D}}

\newcommand{\sT}{\mathcal{T}}

\newcommand{\R}{\mathbb{R}}
\newcommand{\E}{\mathbb{E}}

\newcommand{\bv}[1]{\mathbf{#1}}

\newcommand{\KL}{\text{KL}}

\newcommand{\DKL}{\mathcal{D}_{\KL}}

\newtheorem{Remark}{Remark}

\newtheorem{Definition}{Definition}

\newtheorem{Assumption}{Assumption}

\newtheorem{Problem}{Problem}

\newtheorem*{myexp*}{Running example: control of a pendulum (continue)}
\newtheorem{myexpnum}{Example}

\usepackage{lineno,hyperref}

\journal{Annual Reviews in Control}

\makeatletter
\def\ps@pprintTitle{%
  \let\@oddhead\@empty
  \let\@evenhead\@empty
  \let\@oddfoot\@empty
  \let\@evenfoot\@oddfoot
}
\makeatother

\usepackage[symbol]{footmisc}

\renewcommand{\thefootnote}{\fnsymbol{footnote}}

\begin{document}

\begin{frontmatter}

\title{Probabilistic design of optimal sequential decision-making algorithms in learning and control\let\thefootnote\relax\footnotetext{{\bf This is an authors’ version of the work that is published in Annual Reviews in Control,  Vol. 54, 2022, Pages 81-102.  Changes were made to this version by the publisher prior to publication. The final version of record is available at https://doi.org/10.1016/j.arcontrol.2022.09.003}}}

\author{\'Emiland Garrab\'e, Giovanni Russo\fnref{myfootnote}}
\address{Dept. of Information and Electrical Engineering \& Applied Mathematics, University of Salerno}
\fntext[myfootnote]{Email: \url{{egarrabe,giovarusso}@ unisa.it}}

\begin{abstract}
This survey is focused on certain sequential decision-making problems that involve optimizing over probability functions. We discuss the relevance of these problems for  learning and control.   The survey is organized around a framework that combines a problem formulation and a set of resolution methods.  The formulation consists of an infinite-dimensional optimization problem.  The methods  come from approaches to search optimal solutions in the space of probability functions.  Through the lenses of this overarching framework we revisit popular learning and control algorithms, showing that these naturally arise from suitable variations on the formulation mixed with different resolution methods.  A running example, for which we make the code available, complements the survey.  Finally, a number of  challenges arising from the survey are also outlined.
\end{abstract}


\end{frontmatter}


\section{Introduction}
Sequential decision-making (DM) problems  are ubiquitous in many scientific domains, with application areas spanning e.g., engineering,  economics, management and health \cite{Powell_book,Sharing_book,meyn_2022}. These problems involve a feedback loop where, at each time-step, a decision-maker determines a decision based on the available information.  The result, from the viewpoint of an external observer,  is a sequence of sensed information and decisions that are iterated over time. 

Given their relevance to a wide range of applications, there is then no surprise that, over the years, several communities have worked to address a variety of DM problems, with each community often developing their own toolkit of techniques to tackle the formulations of their interest.  In this survey,  we revisit certain sequential DM problems having probability functions as decision variables and discuss how these problems naturally arise in certain reinforcement learning (RL) and control domains, including the emerging data-driven control (DDC) domain, that have a randomized policy as optimal solution.  The survey is organized around  a framework that consists of a problem formulation and of a set of methods to tackle the problem.  In turn, the formulation consists of an infinite-dimensional optimization problem having probability functions as decision variables.  The methods come from {\em ideas} to search the optimal solution through probability functions.  Equipped with the framework,  we show that popular learning and control algorithms arise from mixing  different variations o{f} the formulation with different resolution methods. The survey is complemented with a tutorial element\footnote{The code to replicate all the numerical results  is made openly available at \url{https://github.com/GIOVRUSSO/Control-Group-Code/tree/master/Decision-making.}}: by developing a running example we illustrate the more applied aspects of certain resolution methods, highlighting some of the key algorithmic details.  The framework, together with the running example, also leads us to highlight a number of application and methodological challenges.  

The paper is organized as follows.  We first (Section \ref{sec:set-up}) introduce the mathematical set-up and formulate the decision-making problem used as an overarching framework within the survey.  Probability functions are central to the proposed formulation and hence, in Section \ref{sec:pdfs}, we expound certain links between these functions and  stochastic/partial differential equations (i.e., SDEs and PDEs). In the same section we also report a conceptual algorithm to compute probabilities from data.  Once this framework is introduced, we survey a set of techniques to solve the problem in the context of learning and control (respectively, in Section \ref{sec:RL} {where multi-armed bandits are also covered} and Section \ref{sec:control}).   We use as a running example {the control of an inverted pendulum} to complement our discussion.  Concluding remarks are given in Section \ref{sec:conclusions}. 

\section{The set-up}\label{sec:set-up}

Vectors are denoted in {\bf bold}.  Let $\mathbb{N}_0$ be the set of positive integers, $\mathcal{X} \subseteq \mathbb{R}^{n_x}, n_x \in \mathbb{N}_0$ and $\mathcal{F}$ be a $\sigma$-algebra on $\mathcal{X}$. Then, the (vector) random variable on $(\mathcal{X},\mathcal{F})$ is denoted $\mathbf{X}$ and its realization by $\mathbf{x}$.   The expectation of a function $f:\mathcal{X}\rightarrow \mathbb{R}$ is denoted by $\mathbb{E}_p\left[f(\mathbf{X})\right]$, where $p(\mathbf{x})$ is the probability density function (if $\mathbf{X}$ is continuous) or probability mass function (if it is discrete) of $\mathbf{X}$. We use the notation $\bv{x} \sim p(\bv{x})$ to state that $\bv{x}$ is sampled from $p(\bv{x})$.  In what follows, we simply say that $p(\bv{x})$ is a {\em probability function} (pf) and we denote by $\mathcal{S}(p)$ its support.  For example,  in what follows $\mathcal{N}(\mu,\sigma)$ denotes a Gaussian (or normal) pf with mean $\mu$ and variance $\sigma$ (the support of the Gaussian is $\R$). The joint pf of two random variables, $\mathbf{X}_1$ and $\mathbf{X}_2$, is written as $p(\mathbf{x}_1, \mathbf{x}_2)$ and the conditional pf of $\mathbf{X}_1$ given $\mathbf{X}_2$ is denoted by $p(\mathbf{x}_1\mid\mathbf{x}_2)$.  Whenever we consider integrals and sums involving pfs we always assume that the integrals/sums exist. Functionals are denoted by capital calligraphic letters with their arguments in curly brackets. In what follows, the convex set of probability functions is denoted by $\mathcal{P}$. We make use of the Matlab-like notation $k_1:k_2$ and $\mathbf{x}_{k_1:k_2}$, with $k_1 \ge k_2$ being two integers, to compactly denote the ordered set of integers $\left\{{k_1},\ldots, {k_2}\right\}$ and the ordered set $\left\{\mathbf{x}_{k_1},\ldots,\mathbf{x}_{k_2}\right\}$, respectively. Following the same notation, we denote by $\left\{p_k(\bv{x}_k)\right\}_{k\in k_1:k_2}$ the ordered set $\left\{p_{k_1}(\bv{x}_{k_1}),\ldots,p_{k_2}(\bv{x}_{k_2})\right\}$.  Subscripts  denote the time-dependence of certain variables; in particular, we use the subscript $k$ for variables that depend on time {discretely} and the subscript $t$ for variables that depend on the time continuously. Finally,  we denote the \textit{indicator function} of $\mathcal{X}$ as $\mathds{1}_{\mathcal{X}}(\mathbf{x})$.  That is, $\mathds{1}_{\mathcal{X}}(\mathbf{x}) =1$, $\forall \mathbf{x} \in \mathcal{X}$ and $0$ otherwise. 

\subsection{Probability functions as a way to describe closed-loop systems}\label{sec:closed_loop_systems}

\begin{figure}[thbp]
\centering
\tikzset{every picture/.style={line width=0.75pt}} 

\begin{tikzpicture}[x=0.75pt,y=0.75pt,yscale=-1,xscale=1]

\draw   (223.8,204) -- (370.8,204) -- (370.8,233.06) -- (223.8,233.06) -- cycle ;
\draw [color={rgb, 255:red, 253; green, 0; blue, 0 }  ,draw opacity=1 ]   (207,222) -- (224,222) ;
\draw [color={rgb, 255:red, 255; green, 0; blue, 0 }  ,draw opacity=1 ]   (207,131) -- (207,222) ;
\draw    (400,132) -- (338,132) ;
\draw    (400,132) -- (400,222) ;
\draw    (400,222) -- (374,222) ;
\draw [shift={(372,222)}, rotate = 360] [color={rgb, 255:red, 0; green, 0; blue, 0 }  ][line width=0.75]    (10.93,-3.29) .. controls (6.95,-1.4) and (3.31,-0.3) .. (0,0) .. controls (3.31,0.3) and (6.95,1.4) .. (10.93,3.29)   ;
\draw [color={rgb, 255:red, 255; green, 0; blue, 0 }  ,draw opacity=1 ]   (207,131) -- (256,131) ;
\draw [shift={(258,131)}, rotate = 180] [color={rgb, 255:red, 255; green, 0; blue, 0 }  ,draw opacity=1 ][line width=0.75]    (10.93,-3.29) .. controls (6.95,-1.4) and (3.31,-0.3) .. (0,0) .. controls (3.31,0.3) and (6.95,1.4) .. (10.93,3.29)   ;
\draw [color={rgb, 255:red, 65; green, 117; blue, 5 }  ,draw opacity=1 ]   (338.8,142) -- (390,142) ;
\draw [color={rgb, 255:red, 65; green, 117; blue, 5 }  ,draw opacity=1 ]   (390,211) -- (390,142) ;
\draw [color={rgb, 255:red, 65; green, 117; blue, 5 }  ,draw opacity=1 ]   (390,211) -- (375,211) ;
\draw [shift={(373,211)}, rotate = 360] [color={rgb, 255:red, 65; green, 117; blue, 5 }  ,draw opacity=1 ][line width=0.75]    (10.93,-3.29) .. controls (6.95,-1.4) and (3.31,-0.3) .. (0,0) .. controls (3.31,0.3) and (6.95,1.4) .. (10.93,3.29)   ;
\draw   (260,118) -- (337.8,118) -- (337.8,147.06) -- (260,147.06) -- cycle ;

\draw (267,121) node [anchor=north west][inner sep=0.75pt]   [align=left] {{\large System}};
\draw (233,208) node [anchor=north west][inner sep=0.75pt]   [align=left] {{\large Decision-maker}};
\draw (177,162) node [anchor=north west][inner sep=0.75pt]    {$\textcolor[rgb]{1,0,0}{\mathbf{u}}\textcolor[rgb]{1,0,0}{_{k}}$};
\draw (403,163) node [anchor=north west][inner sep=0.75pt]    {$\mathbf{x}_{k-1}$};
\draw (244,86) node [anchor=north west][inner sep=0.75pt]    {$f_{\mathbf{x} ,k}(\mathbf{x}_{k} \mid \mathbf{x}_{_{k-1}} ,\textcolor[rgb]{1,0,0}{\mathbf{u}}\textcolor[rgb]{1,0,0}{_{k}})$};
\draw (249,245) node [anchor=north west][inner sep=0.75pt]    {$\textcolor[rgb]{0,0,0}{f_{\mathbf{u} ,k}(}\textcolor[rgb]{1,0,0}{\mathbf{u}}\textcolor[rgb]{1,0,0}{_{k}} \mid \mathbf{x}_{k-1})$};
\draw (334,162) node [anchor=north west][inner sep=0.75pt]    {$\textcolor[rgb]{0.25,0.46,0.02}{r}\textcolor[rgb]{0.25,0.46,0.02}{_{k}}\textcolor[rgb]{0.25,0.46,0.02}{\ or\ c}\textcolor[rgb]{0.25,0.46,0.02}{_{k}}$};

\end{tikzpicture}

\caption{the decision-making feedback loop. The decision at time-step $k$, $\bv{u}_k$, is determined {from information available at $k-1$; $\bv{x}_{k-1}$ denotes the state at $k-1$}.  Decisions are driven by $r_k$ ($c_k$), i.e., the reward (cost) received at each $k$ (see Section \ref{sec:gen_prob_statement}).  The pfs from which $\bv{u}_k$ and $\bv{x}_{k-1}$ are sampled are introduced in Section \ref{sec:closed_loop_systems}.
}
\label{fig:closed_loop}
\end{figure}
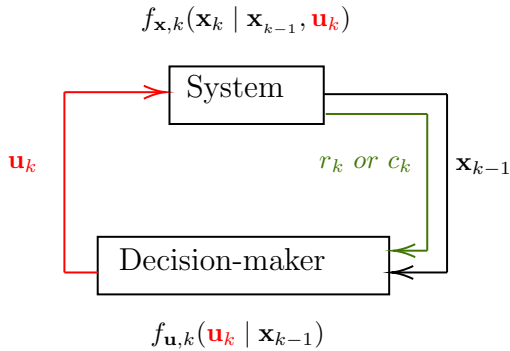

We consider the feedback loop (or closed-loop system) schematically illustrated in Figure \ref{fig:closed_loop}, where a {\em decision-maker} interacts with the {\em system} with the goal of fulfilling a given task. In certain applications closer to the RL community, the decision-maker is termed as an agent and the system with which it interacts is the environment \cite{RL_intro_2018}.  Within the control community, the decision-maker is typically a control algorithm and the system is often the plant under control. {The terms control inputs/actions/decisions and agent/decision-maker/control algorithm are used interchangeably throughout this paper}.

Consider the time-horizon $\sT:= 0:T$ and let: (i) $\bv{U}_k$ be the action determined by the decision-maker at time-step $k$; (ii) $\bv{X}_{k-1}$ be the state variable at time-step $k-1$.  {Fundamental to properly formalize, and study, DM problems is the concept of state}. This variable represents what the decision-maker knows in order to determine $\bv{U}_k$. {In certain streams of literature on DM surveyed in \cite[Chapter $9$]{Powell_book}, the state is explicitly partitioned as $\bv{X}_k := (\bv{O}_k,\bv{B}_k)$,  where:
\begin{itemize}
    \item $\bv{O}_k$ embeds all observable states (i.e., the explicitly known information) that are needed to make a decision;
    \item $\bv{B}_k$ is the belief and specifies pfs that describe quantities that are not directly observed. These, while not explicitly known, are needed to make a decision (see also our discussion on partial observations in Section \ref{sec:POMDPs}).
\end{itemize}}
Both the state and the action can be vectors and we let $n_x$ and $n_u$ be the dimensions of these vectors.  {For our derivations in what follows it is useful to introduce the following dataset \begin{equation}\label{eqn:dataset}
\boldsymbol{\Delta}_{0:T}{:=\{\mathbf{x}_0, \mathbf{u}_1, \mathbf{x}_1, \mathbf{u}_2, \ldots, \mathbf{x}_{T-1}, \mathbf{u}_T, \mathbf{x}_{T}\}},
\end{equation}}
collected from the closed-loop system of Figure \ref{fig:closed_loop} over $\sT$. {We note that, while the rewards/costs received by the agent at each $k$ do not explicitly appear in \eqref{eqn:dataset}, {these \emph{drive} the DM process. The reward/cost received by the agent at each $k$ influences the decisions of the decision-maker and hence the future evolution of the state.} This dependency can be explicitly highlighted by introducing an {\em exogenous} random variable, which can be used to capture uncertain information received by the decision-maker from the system/environment (hence including the reward/cost signal received by the agent). In particular, let $\bv{W}_k$ be this exogenous random variable at time-step $k$, sometimes in the literature the notation $\bv{x}_k(\bv{w}_k)$ is used to stress that the state depends on this exogenous information. Similarly, one can write $\bv{u}_k(\bv{w}_{k-1})$ to stress that the decision made by the decision-maker depends on the exogenous random variable.}  
Then, the evolution of the closed-loop system can be described \cite{bayesian_id,TOUCHETTE2004140}  by the joint pf, say $f(\boldsymbol{\Delta}_{0:T}).$ By making the standard Markov assumption and by applying the chain rule for pfs \cite{bayesian_id} we obtain the following convenient factorization for $f(\boldsymbol{\Delta}_{0:T})$:
\begin{equation}\label{eqn:full_joint_general}
    f(\boldsymbol{\Delta}_{0:T}) = f_{0}(\mathbf{x}_0)\prod_{k\in1:T}f_{\bv{x},k}(\mathbf{x}_{k}\mid\mathbf{x}_{k-1}, \mathbf{u}_{k})f_{\bv{u},k}(\mathbf{u}_{k}\mid\mathbf{x}_{k-1}).
\end{equation}
We refer to (\ref{eqn:full_joint_general}) as the probabilistic description of the closed loop system.  The pfs $f_{\bv{x},k}(\mathbf{x}_{k}\mid\mathbf{x}_{k-1}, \mathbf{u}_{k})$ describe {how the state evolves at each $k$.} This is termed as the probabilistic description of the system and we denote its support by $\mathcal{X}\subseteq\R^{n_x}$. The control input at time-step $k$ is determined by sampling from $f_{\bv{u},k}(\mathbf{u}_{k}\mid \mathbf{x}_{k-1})$. This is a randomized policy: {it is the probability of making decision $\bv{u}_k$ given $\bv{x}_{k-1}$. In what follows,} $f_{\bv{u},k}(\mathbf{u}_{k}\mid \mathbf{x}_{k-1})$ is termed as control pf and has support $\mathcal{U}\subseteq\R^{n_u}$.  In (\ref{eqn:full_joint_general}) the initial conditions are embedded via  the prior $f_{0}(\mathbf{x}_0)$.  For time invariant systems $f_{\bv{x},k}(\mathbf{x}_{k}\mid\mathbf{x}_{k-1}, \mathbf{u}_{k})$ in (\ref{eqn:full_joint_general}) is the same $\forall k$.  Also, a policy is  stationary if $f_{\bv{u},k}(\mathbf{u}_{k}\mid \mathbf{x}_{k-1})$ is the same $\forall k$. In both cases,  when this happens we drop the $k$ from the subscripts in the notation. 

\begin{Remark}
The time-indexing used in (\ref{eqn:full_joint_general}) is chosen in a way such that, at each $k$, the decision-maker determines $\bv{u}_k$ based on data at $k-1$. With this time indexing, once the system receives $\bv{u}_k$, its state transitions from $\bv{x}_{k-1}$ to $\bv{x}_{k}$. {{When} the exogenous random variables are used in the notation, the time indexing is such that $\bv{W}_{k}$ is available to the agent when action $\bv{u}_k$ is determined.} 
\end{Remark}
{\begin{Remark}\label{rem:Wk}
The pfs formalism in \eqref{eqn:full_joint_general} leveraged in this paper is typically used within the literature on Markov Decision Processes (MDPs) and sequential DM under uncertainty. The pf $f_{\bv{u},k}(\mathbf{u}_{k}\mid \mathbf{x}_{k-1})$ is the probability of making decision $\bv{u}_k$ given $\bv{x}_{k-1}$. This probability depends on the exogenous information, i.e., the $\bv{W}_k$'s, received by the decision-maker. This includes the reward/cost received, at each $k$, by the agent. The pf $f_{\bv{x},k}(\mathbf{x}_{k}\mid\mathbf{x}_{k-1}, \mathbf{u}_{k})$ capturing the evolution of $\bv{X}_k$ can be expressed via a dynamics of the form $\bv{X}_k = f_k(\bv{X}_{k-1},\bv{U}_k,\bv{W}_k)$. In e.g., \cite[Chapter $9$]{Powell_book} it is shown that pfs can be computed from the dynamics (see also Section \ref{sec:pdfs} where we discuss the link between pfs and differential equations). 
\end{Remark}}
\begin{Remark}
As remarked in \cite{bayesian_id}, the pf in (\ref{eqn:full_joint_general}) is the most general description of a system from the viewpoint of an outer observer.  In Section \ref{sec:pdfs} we discuss what these pfs can capture and report an algorithm to estimate these pfs from data.  
\end{Remark}
{\begin{Remark}\label{rem:Markov}
Making the Markov assumption is often a {\em design} choice to simplify the resolution of the DM problem. If the Markov property is not satisfied, say $\bv{U}_k$ depends on the past history starting from some $k-\tau$, $\tau >1$, then one can redefine the state to include, at each $k$, all past history up to  $k-\tau$. This choice is often avoided in practice as the idea is to restrict the policy design to some sufficient statistic that does not require storing long histories to compute $\bv{U}_k$. This aspect is related to the notion of information state. In e.g., \cite{AIS}, this is defined as a compression of the system history that is sufficient to predict the reward and the system next state once an action is taken. That is, intuitively, an information state is a statistics that is sufficient for performance evaluation. 
\end{Remark}}
\begin{Remark}
We use the wording {\em dataset} to denote a sequence of data.  In certain applications, where{,} e.g.{,}  multiple experiments can be performed,  one might have access to a collection of datasets. This is termed as {\em database}.   Finally, we use the wording {\em data-point} {to denote the data collected at a given $k$}.
\end{Remark}

\subsection{Statement of the decision-making problem}\label{sec:gen_prob_statement}
The following finite-horizon, infinite-dimensional, sequential decision-making problem serves as an overarching framework for this survey:

\begin{Problem}\label{prob:general_problem}
Let, $\forall k\in 1:T$: 
\begin{enumerate}
    \item $\mathcal{E}_k$ and $\mathcal{I}_k$ be index sets {at time-step $k$};
    \item $H_{\bv{u},k}^{(i)}$,  $G_{\bv{u},k}^{(j)}$, $0\le\varepsilon_k\le 1$ with ${i}\in\mathcal{E}_k$ and ${j}\in\mathcal{I}_k$, be  constants;
    \item $h_{\bv{u},k}^{(i)}$,  $g_{\bv{u},k}^{(j)}:\mathcal{U}\rightarrow \R$, with ${i}\in\mathcal{E}_k$ and ${j}\in\mathcal{I}_k$, be  measurable mappings;
    \item $\bar{\mathcal{X}}_k\subseteq\mathcal{X}$.
\end{enumerate}  
Find  $\{f^{\ast}_{\bv{u},k}(\mathbf{u}_{k}\mid \mathbf{x}_{k-1})\}_{k \in 1:T}$ such that:
\begin{subequations}
    \begin{alignat}{3}
    \{f_{\bv{u},k}^{\ast}(\mathbf{u}_{k}\mid \mathbf{x}_{k-1})\}_{k \in 1:T} & \in \underset{\{f_{\bv{u},k}(\mathbf{u}_{k}\mid \mathbf{x}_{k-1})\}_{k \in 1:T}}{{\arg \min}}\  \mathbb{E}_{f}[c_{{1:T}}(&&\mathbf{X}_{0},\ldots,\bv{X}_T, \mathbf{U}_{1}, \ldots,\bv{U}_T)]\label{genProb_cost}\\
    & s.t.   \ \mathbf{x}_{k} \sim f_{\bv{x},k}(\mathbf{x}_{k}\mid \mathbf{x}_{k-1}, \mathbf{u}_{k}),   &&\forall k\in 1:T;\label{genProb_xSamp}\\
    &   \ \ \ \ \ \mathbf{u}_{k} \sim f_{\bv{u},k}(\mathbf{u}_{k}\mid \mathbf{x}_{k-1}), &&\forall k\in 1:T;\label{genProb_uSamp}\\
   &   \ \ \ \ \  \E_{f_{\bv{u},k}}\left[h^{({i})}_{\bv{u},k}(\bv{U}_k)\right] = H_{\bv{u},k}^{{(i)}}, &&\forall k\in 1:T,    {\forall i}\in \mathcal{E}_k;\label{genProb_eqConstr}\\
    &   \ \ \ \ \  \E_{f_{\bv{u},k}}\left[g^{(j)}_{\bv{u},k}(\bv{U}_k)\right] \le G_{\bv{u},k}^{{(j)}}, &&\forall k\in 1:T,   \forall j\in \mathcal{I}_k;\label{genProb_ineqConstr}\\
    &  \ \ \ \ \  \mathbb{P}\left(\bv{X}_k\in\bar{\mathcal{X}}_k\right) \ge 1 -\varepsilon_k, &&\forall k\in 1:T;\label{genProb_safetySet}\\
        &   \ \ \ \ \ f_{\bv{u},k}(\mathbf{u}_{k}\mid \mathbf{x}_{k-1}) \in \mathcal{P}, &&\forall k\in 1:T.\label{genProb_pf}
    \end{alignat}
\end{subequations}
\end{Problem}
In the cost of Problem \ref{prob:general_problem} the expectation is over the pf $f(\boldsymbol{\Delta}_{0:T})$ {and in the problem statement we used the shorthand notation $f$ to denote this pf}. A typical assumption is that, at each $k$, the cost depends only on $\bv{X}_{k-1}$ and $\bv{U}_k$ (the implications of this assumption are exploited in Section \ref{sec:RL} and Section \ref{sec:control}). As usual, we express the DM problem via the minimization of a cost (or maximizing a reward). The assumption that {\em any} DM problem can be stated via some reward maximization is known as the {\em reward hypothesis}. We refer readers to e.g., \cite{Deb_54,SILVER2021103535} for a detailed discussion on the validity of such hypothesis.
\begin{Remark}
Intuitively,  the minimization over policies in Problem \ref{prob:general_problem} means that the goal of the optimization problem is that of finding the {\em best method} for making decisions.  As we shall see,  certain decision-making approaches rely on optimizing directly over the control variable. In turn, this means finding the {\em best action} to optimize the cost and not the method that generated the action.
\end{Remark}

In Problem \ref{prob:general_problem} the decision variables are the pfs  $f_{\bv{u},k}(\mathbf{u}_{k}\mid \mathbf{x}_{k-1})$. Constraints (\ref{genProb_xSamp}) - (\ref{genProb_uSamp}) capture the fact that, at each $k$, the  state and control are sampled from the probabilistic description of the system and the control pf (as we shall see, the pfs $f_{\bv{x},k}(\mathbf{x}_{k}\mid \mathbf{x}_{k-1}, \mathbf{u}_{k})$ do not necessarily need to be known).  {The sets $\mathcal{E}_k$ and $\mathcal{I}_k$ in the problem statement are the index sets at time-step $k$ for the equality and inequality actuation constraints (\ref{genProb_eqConstr}) - (\ref{genProb_ineqConstr}). That is, we used the superscripts $(i)$ to denote constants and mappings related to the $i$-th equality constraint (analogously, the superscript $(j)$ is related to the $j$-th inequality constraint). In the formulation, the actuation constraints are formalized as expectations of the (possibly nonlinear) mappings $h_{\bv{u},k}^{(i)}$ and $g_{\bv{u},k}^{(j)}$. Hence, even if the mappings can be nonlinear, the constraints are linear in the decision variables of Problem \ref{prob:general_problem}. Also, $H_{\bv{u},k}^{(i)}$ and  $G_{\bv{u},k}^{(j)}$ are constants appearing on the right hand side of (\ref{genProb_eqConstr}) - (\ref{genProb_ineqConstr}).  Note that the formulation allows to consider situations where  the equality and inequality constraints, and their number, can change over time (see items $1$ - $3$ at the beginning of Problem \ref{prob:general_problem}). The constraints (\ref{genProb_eqConstr}) - (\ref{genProb_ineqConstr})} can be used to guarantee properties on the moments of the decision variables {(see \cite{gagliardi2020probabilistic} for a detailed discussion)}. Constraint (\ref{genProb_ineqConstr}) can also capture bound constraints on $\bv{U}_k$ of the form $\mathbb{P}(\bv{U}_k \in \bar{\mathcal{U}}) \ge 1 - \gamma_k$, where $\bar{\mathcal{U}} \subseteq \mathcal{U}$ and $0\le\gamma_k \le 1$. See \cite{gagliardi2020probabilistic} for a detailed discussion, where it is also shown that these constraints are convex in the decision variables and hence can be solved without resorting to bounding approximations. The fulfillment of constraint (\ref{genProb_safetySet}) instead guarantees that the probability that the state is outside some (e.g., desired) set ${\bar{\mathcal{X}}}_k$ is {less} than some {\em acceptable} $\varepsilon_k$. Finally, the last constraint is a {\em normalization} constraint guaranteeing that the solution to the problem belongs to {the convex set of pfs} $\mathcal{P}$ {(see also our notation at the beginning of Section \ref{sec:set-up})}.

Problem \ref{prob:general_problem} allows to consider problems with both continuous and discrete states/controls. In the former case, the pfs are probability density functions, while in the latter these are probability mass functions.  Interestingly, the  formulation also allows to consider situations where the available control inputs are {\em symbols}, e.g., {\em go left}, {\em go right}.  Indeed, in this case  $f_{\bv{u},k}(\mathbf{u}_{k}\mid \mathbf{x}_{k-1})$ can be defined over the list of possible symbols.  See also \cite{Todorov_pnas} for a related discussion.  Rather counter-intuitively,  as we discuss in Section \ref{sec:control},  analytical expressions can be obtained for the optimal solution of relaxed versions of Problem \ref{prob:general_problem} even when the actions are symbolic. Finally, in the context of neuroscience \cite{Pir_21}  a special case of Problem \ref{prob:general_problem} (i.e., without the second and third constraints and with a cost split into short-term and long-term components) is used to suggest a solution for a central puzzle about the neural representation of cognitive maps in humans. 

\begin{Remark}
In our problem formulation  we optimize over pfs. Formally, deterministic policies in the state variable can be written as pfs.  Hence, in principle,  Problem \ref{prob:general_problem} can be used to study both randomized and deterministic policies. Nevertheless,  here we also consider control problems that have as optimal solution a randomized policy (that is, a pf). These problems typically go under the label of probabilistic control/design problems and  are discussed in Section \ref{sec:control}.
\end{Remark}


{
\subsection{A discussion on partial observations}\label{sec:POMDPs}
Situations where only partial {observations} (rather than the full state) are available to the decision-maker naturally arise in a number of applications, such as robotic planning \cite{sampling_hyperbelief}, finance and healthcare \cite{AIS}. 
In this context, we note that a number of techniques are available to reduce DM problems with partial observations to (fully observed) MDPs whose state is a belief state. The belief at time-step $k$ describes the distribution of the state given the information available to the decision-maker up to $k$ \cite{pomdp_to_mdp}. Note that, as discussed in Section \ref{sec:closed_loop_systems}, the presence of belief state does not preclude the existence of states that can be directly observed by the agent so that a subset of the state variables is directly observable, while other state variables are represented via their belief. We also highlight that, in certain streams of the literature, belief states are leveraged to represent some {parameter} of the system that is unknown. We refer to \cite[Chapter $9$]{Powell_book} for a discussion on this aspect -- in such a work it is also discussed how solving DM problems with partial observations and large belief states can become intractable using classic resolution approaches for the fully observed set-up. As noted in \cite{perseus_pomdp}, in partially observable environments, some form of memory is needed in order for the agent to compute its decisions. As also discussed in this work, if the transition and observation models of a partially observed MDP are known, then this can be recast into a belief-state MDP. In this context, in \cite{perseus_pomdp} a randomized point-based value iteration algorithm, {\em PERSEUS}, for partially observed MDPs is introduced. PERSEUS can operate on a large belief space and relies on simulating random interactions of the agent with the environment. Within the algorithm, a number of value backup stages are performed and it is shown that  in each backup stage the value of each point in the belief space does not decrease. A complementary line of research, inspired by graph-sampling methods, can be found in \cite{sampling_hyperbelief}. In this work, optimal plans of decisions are searched in the hyperbelief space (i.e., the space of pfs over the belief) using an approach devised from these methods. In particular, the problem is abstracted into a two-level {\em planner} and the approach, which leverages a graph representation in the hyperbelief space, is shown to have the following features: (i) optimization over the graph can be performed without exponential explosion in the dimension of the hyperbelief space; (ii) the bound on the optimal value can only decrease at every iteration. When the pf describing the evolution of the system/environment is not available, techniques known as partially-observed RL have been developed. For example, in \cite{AIS} it is shown that if a function of the {\em history} approximately satisfies the properties of the information state (see also Remark \ref{rem:Markov}) then the policy can be computed using an approximate dynamic programming decomposition. Moreover, the policy is approximately optimal with bounded loss of optimality. Finally, we recall that in \cite{po_rl_scary} it is shown that for a wide class of partially-observed RL problems, termed as {weakly revealing partially observed MDPs}, learning can be made sample-efficient.

}

\section{Relating probability functions to SDEs, PDEs and data}\label{sec:pdfs}
Probability functions are central to the formulation of  Problem \ref{prob:general_problem} and now we briefly expound certain links between pfs,  SDEs and PDEs.  We also report a conceptual algorithm to estimate pfs from data. We start with  considering the SDE in the It\= o sense (satisfying the usual conditions on the existence and uniqueness of the solutions) of the form:
\begin{equation}\label{eqn:Ito_SDE}
d\bv{X}_t = b(\bv{X}_t,t)dt + \sigma(\bv{X}_t,t)d\bv{W}_t,
\end{equation}
where  $\bv{X}_t\in\mathcal{X}\subseteq\R^{n_x}$,  $\bv{W}_t$ is an $n_w$-dimensional Wiener process, $b(\cdot,\cdot)$ is the drift function and $\sigma(\cdot,\cdot)$ is the $n_x\times n_w$ full-rank diffusion matrix.  The solution of (\ref{eqn:Ito_SDE}) is a Markov process  (see e.g., Theorem $9.1$ in \cite{Mao_97}) characterized by the transition density probability function $\hat{\rho}(\bv{x},t;\bv{y},s)$. This is the transition density probability function for the stochastic process to move from $\bv{y}$ at time $s$ to $\bv{x}$ at time $t$. The Fokker-Planck (FP) equation \cite{fokker,planck} associated to (\ref{eqn:Ito_SDE}) is given by:
\begin{equation}\label{eqn:FP_general}
\partial_t \rho(\bv{x},t) + \sum_{i{\in 1:n_x}} \partial_{x_i}\left(b_i(\bv{x},t)\rho(\bv{x},t)\right) - \sum_{i,j{\in 1:n_x}}\partial^2_{x_ix_j}\left(a_{ij}(\bv{x},t)\rho(\bv{x},t)\right) = 0,
\end{equation}
where {$\rho(\bv{x},t)$} is the probability density to find the process (\ref{eqn:Ito_SDE}) at $\bv{x}$ at time $t$.  In the above expression the subscripts denote the elements of vectors/matrices.  The FP equation is a PDE of parabolic type with its Cauchy data given by the initial pf $\rho(\bv{x},0) = \rho_0(\bv{x})$. The diffusion coefficients in (\ref{eqn:FP_general}) are $a_{ij}(\bv{x},t) := \frac{1}{2} \sum_{k{\in 1:n_x}}\sigma_{ik}(\bv{x},t)\sigma_{jk}(\bv{x},t)$. In  Section \ref{sec:FP_control} we survey a set of methods that exploit the link between pfs, SDEs and the FP equation. 

\begin{Remark}
Besides capturing physical processes governed by PDEs and SDEs,  in e.g., \cite{Todorov_pnas} it is noted how pfs can be leveraged to capture the evolution of processes that have discrete and/or symbolic states and inputs (see also our discussion at the end of Section \ref{sec:set-up}). Further, the pfs formalism also naturally arises in probabilistic programming, as well as in applications where a given system can be modeled via probabilistic Boolean networks or, in a broader context, via Markov random fields and Bayesian networks, see e.g., \citep{10.1145/2593882.2593900,10.1093/bioinformatics/18.2.261,10.5555/1051482,10.5555/1535531}.
\end{Remark}

In a broad sense, the problem of estimating pfs that fit a given set of data goes under the label of {\em density estimation}.  While surveying methods to estimate densities goes beyond the scope of this paper, we refer readers to \cite{Silverman:1070306} for a detailed survey of different techniques and to \cite{10.5555/1121596} for applications to robotics.  For completeness, we also report an algorithm to estimate conditional pfs from data that is used within our illustrative examples.  The pseudo-code for the algorithm, which is adapted from   histogram filters, is given in Algorithm \ref{alg:hist_filter}. The algorithm is a non-parametric filter to estimate the generic pf $p(\mathbf{z}_k\mid\mathbf{y}_{k-1})$ from a sequence of data $\{(\bv{z}_k,\bv{y}_{k-1})\}_{1:N}$,  where $\bv{z}_k\in\mathcal{Z}$ and $\bv{y}_{k-1}\in\mathcal{Y}$. This is done by first discretizing $\mathcal{Z}$ and $\mathcal{Y}$ (steps $3-4$) and then by binning the data to obtain the empirical joint pfs $p(\bv{y}_{k-1})$ and $p(\bv{z}_k,\bv{y}_{k-1})$. This latter operation is done in steps $5-6$ (note that the binning function provides a normalized histogram and takes as input both the data and the discretized sets over which the binning is performed). Once the joint pfs are computed,  the estimate is obtained via Bayes rule. This is done in steps $7-15$, where: (i)  a logical condition is included, which sets $p(\hat{\bv{z}}_k\mid\hat{\bv{y}}_{k-1})$ to $0$ whenever $p(\hat{\bv{y}}_{k-1})=0$, i.e., whenever the event $\bv{Y}_{k-1} = \hat{\bv{y}}_{k-1}$ is not contained in the data; (ii) it is implicit, in step $12$, that a normalization operation is performed.  Algorithm \ref{alg:hist_filter} is a Bayes filter applied on the binned data and it is interesting to notice that: (i) in the ideal situation where the bin width is $0$, the two algorithms coincide; (ii) popular parametric filters such as Gaussian filters are derived from Bayes filters.

\begin{algorithm}
    \caption{{Histogram filter}}\label{pmf_from_data_binning}
\begin{algorithmic}[1]
\State \textbf{Input:} $\{(\bv{z}_k,\bv{y}_{k-1})\}_{1:N}$, $\bv{z}_k\in\mathcal{Z}$, $\bv{y}_{k-1}\in\mathcal{Y}$
\State \textbf{Output:} An estimate of $p(\mathbf{z}_k\mid\mathbf{y}_{k-1})$
\State $\mathcal{Z}_d\gets\text{discretize}(\mathcal{Z})$
\State $\mathcal{Y}_d\gets\text{discretize}(\mathcal{Y})$
\State $p(\bv{y}_{k-1}) \gets \text{bin}(\bv{y}_{0:N-1},\mathcal{Y}_d)$
\State $p(\bv{z}_k,\bv{y}_{k-1}) \gets \text{bin}(\bv{z}_{1:N},\mathcal{Z}_d,\bv{y}_{0:N-1},\mathcal{Y}_d)$
\For{$\hat{\bv{y}}_{k-1}$ in $\mathcal{Y}_d$}
    \For{$\hat{\bv{z}}_{k}$ in $\mathcal{Z}_d$}
            \If{$p(\hat{\bv{y}}_{k-1}) == 0$}
            \State $p(\hat{\bv{z}}_k\mid\hat{\bv{y}}_{k-1}) \gets 0$
            \Else
            \State $p(\hat{\bv{z}}_k\mid\hat{\bv{y}}_{k-1}) \gets \frac{p(\hat{\bv{z}}_k,\hat{\bv{y}}_{k-1})}{p(\hat{\bv{y}}_{k-1})}$
            \EndIf
  \EndFor
\EndFor
\end{algorithmic}\label{alg:hist_filter}
\end{algorithm}

\begin{figure}[thbp]
\centering
\includegraphics[width=\linewidth]{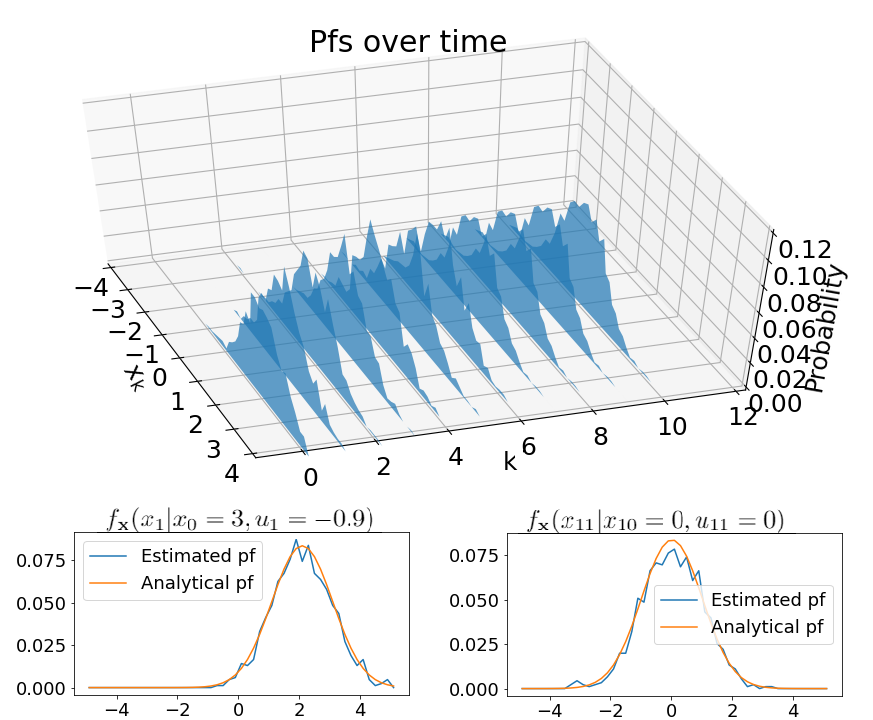}
\caption{evolution of $f_x(x_k\mid x_{k-1},u_{k})$ estimated via Algorithm \ref{alg:hist_filter}.  The figure was obtained by setting the initial condition $X_0=3$. At each $k$, the next state was determined by sampling from $f_x(x_k\mid x_{k-1},0.3x_{k-1})$.  In the bottom panels the estimated pfs at the first and the last time-steps are overlapped to the analytical ones.}
\label{fig:example_comparison_1}
\end{figure}

\begin{figure}[thbp]
\centering
\includegraphics[width=\linewidth]{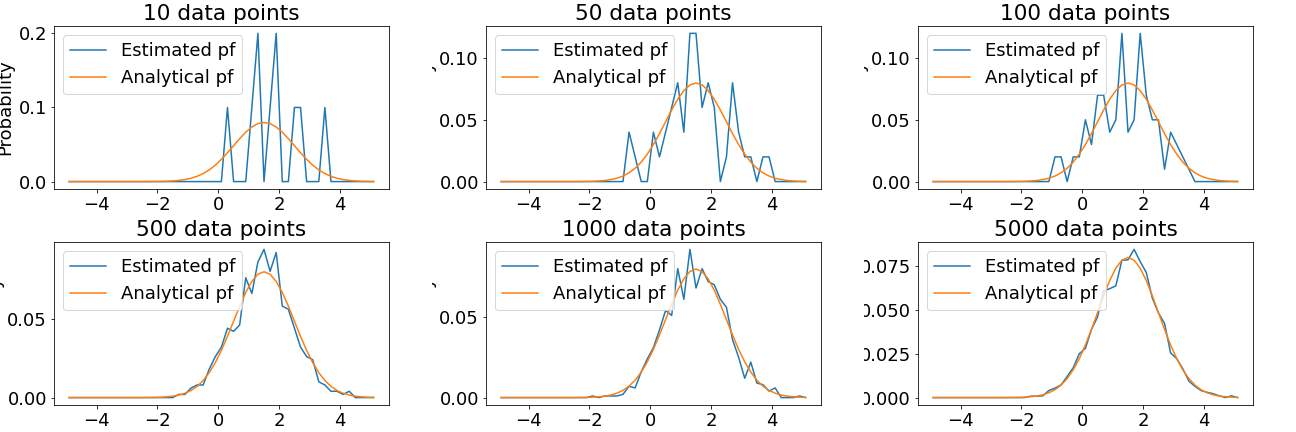}
\caption{illustration of how the estimate via Algorithm \ref{alg:hist_filter} changes as the number of data-points increases.  In each panel, the pf $f_x(x_k\mid1,0.5)$ is shown together with the analytical pf.  The results shown in the panels are representative for all the other pfs. }
\label{fig:side_by_side}
\end{figure}

\begin{myexpnum}\label{exe:pfs}\normalfont

We give a first  example to illustrate how Algorithm \ref{alg:hist_filter} can be used to estimate the pf capturing the evolution of a linear system.  We do so by considering the simple scalar linear system 
\begin{equation}\label{eqn:linear_example}
X_k = X_{k-1} + U_k + W_k,
\end{equation}
with $W_k \sim \mathcal{N}(0,1)$.  It is well known that for the solutions of the above dynamics  it holds that $x_k\sim\mathcal{N}\left(x_{k-1} + u_k,1\right)$ and we now use Algorithm \ref{alg:hist_filter} to estimate this pf from data generated by simulating (\ref{eqn:linear_example}).  To this aim,  we built a database by performing $1000$ simulations (of $100$ time-steps each) of the dynamics in (\ref{eqn:linear_example}). Within each simulation, initial conditions were chosen randomly in the set $[-5,5]$ and, at each time-step of the simulation,  the input was drawn randomly in the set $[-1,1]$. Data-points were removed whenever the state at time-step $k$ fell outside  the range $[-5,5]$.  We discretized both the set $[-5,5]$ for the  state (discretization step of $0.2$) and the range of the inputs $[-1,1]$, with discretization step of $0.1$. In the filter, we also set $\bv{z}_k$ as $x_k$ and $\bv{y}_{k-1}$ as $(x_{k-1}, u_k)$.  This allowed to obtain an estimate of the pf $f_x(x_k\mid x_{k-1},u_k)$. In Figure \ref{fig:example_comparison_1},  a comparison is shown between the estimated pf and the {\em analytical} pf $\mathcal{N}\left(x_{k-1} + u_k,1\right)$ from which $x_k$ is sampled at each time-step.  The figure  illustrates the evolution of the pfs when $u_k$ is generated via the feedback law $U_k = -0.3X_{k-1}$. Finally,  we also numerically investigated how the pfs estimated via Algorithm \ref{alg:hist_filter} change as the number of available data-points increases.  This is reported in Figure \ref{fig:side_by_side}.

\end{myexpnum} 

\section{Running example: control of a pendulum}
We consider the problem of stabilizing a pendulum onto its unstable upward equilibrium.  In this first part of the running example we describe the set-up and the process we followed to compute the pfs.  Our database was obtained by simulating the following discretized pendulum dynamics:
\begin{equation}\label{eqn:pendulum_dynamics}
\begin{split}
\theta_k & = { \theta_{k-1} + } \omega_{k-1}dt + W_{\theta} \\
\omega_k  & =  { \omega_{k-1} + } \left(\frac{g}{l}\sin(\theta_{k-1}) + \frac{U_k}{ml^2}\right)dt + W_{\omega},
\end{split}
\end{equation}
where $\theta_k$ is the angular position of the mass at time-step $k$, $\omega_k$ is the angular velocity, $U_k$ is the torque applied to the hinged end. Also, $W_{\theta}\sim \mathcal{N}\left(0,\sigma_{\theta}\right)$ and $W_{\omega}\sim\mathcal{N}(0,\sigma_{\omega})$, with $\sigma_\theta $ being the variance of the noise on $\theta_k$ and $\sigma_{\omega}$  being the variance of the noise on $\omega_k$. The upward  equilibrium  corresponds to an angular position of $0$.  In the above expression,  $l$ is the length of the rod, $m$ is the weight of the mass,  $g$ is the gravity, $dt$ is the discretization step.  Further, in what follows we set $\mathcal{X} := [-\pi,\pi] \times [-5,5]$ and $\mathcal{U} := [-2.5,2.5]$.  

The pendulum we want to control has parameters $dt = 0.1$s, $l=0.6$m,  $m=1$kg, $\sigma_\theta = 6\pi/180$ (i.e., $3$ degrees) and $\sigma_{\omega} = 0.1$.  We let $\bv{X}_k:=(\theta_k,\omega_k)$ and,  as a first step we simulated the dynamics to obtain a database.  Specifically, the database consisted of data-points collected from $10000$ simulations of the  dynamics. Each simulation consisted of $99$ time-steps (i.e.,  $10$ seconds): initial conditions were randomly chosen and,  at each $k$,  a random input from $\mathcal{U}$ was applied.  The next step was  to estimate $f_{\bv{x}}(\bv{x}_k\mid\bv{x}_{k-1}, u_k)$ and this was done by means of Algorithm \ref{alg:hist_filter} following the process we illustrated in Example \ref{exe:pfs}.  In order to use the algorithm we: (i) set   $\bv{y}_{k-1} := (\bv{x}_{k-1}, u_k)$, $\bv{z}_k := \bv{x}_k$; (ii) discretized $\mathcal{X}$ in a grid of $50\times 50$ bins and $\mathcal{U}$ in $20$ bins (the bin width was uniform).   

\begin{figure}[thbp]
\centering
\includegraphics[width=0.48\linewidth]{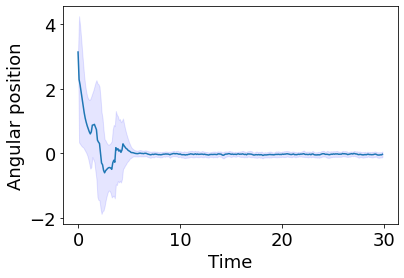}
\includegraphics[width=0.48\linewidth]{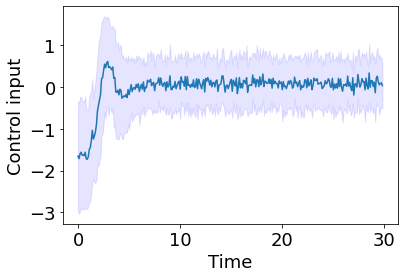}
\caption{behavior of the pendulum obtained by recursively applying $u_k\sim g_{{u}}(u_k\mid\bv{x}_{k-1})$ to (\ref{eqn:pendulum_dynamics}).  Bold lines denote means and shaded areas represent the confidence intervals corresponding to the standard deviation. Figure obtained from $50$ simulations.}
\label{fig:examples_simulation}
\end{figure}

For reasons that will be clear later,  we also obtained a pf for a pendulum that differs from the one considered above in the mass (now weighting $0.5$kg).  When we discuss certain probabilistic methods in Section \ref{sec:control}, this second pendulum will serve as a reference system of which we want to track the evolution.  For the reference system we obtain not only $g_{\bv{x}}(\bv{x}_k\mid\bv{x}_{k-1}, u_k)$ but also a randomized policy, i.e., $g_{{u}}(u_k\mid\bv{x}_{k-1})$,  able to perform the swing-up.  The pf $g_{\bv{x}}(\bv{x}_k\mid\bv{x}_{k-1}, u_k)$ was obtained by following the same process described above for $f_{\bv{x}}(\bv{x}_k\mid\bv{x}_{k-1}, u_k)$. The pf $g_{{u}}(u_k\mid\bv{x}_{k-1})$ was instead obtained by leveraging Model Predictive Control (MPC).  In the MPC formulation: (i) we used the discretized pendulum as model; (ii) the width of the receding horizon window,  $H$, was $20$ time-steps; (iii) at each $k$ the cost within the receding horizon window was $\sum_{t{\in k:k+H-1}}\left(\theta_t^2 + 0.1\omega_t^2\right) + \theta_{t}^2 + 0.5\omega_{k+H}^2$; (iv) the control variable was constrained to belong to the set $[-2,2]$.  Once we obtained a policy from MPC, we added to the control signal a noise process sampled from $\mathcal{N}(0,\sigma_u)$, with $\sigma_u = 0.2$. By doing so, we obtained $g_{{u}}(u_k\mid\bv{x}_{k-1})$: by construction this is a Gaussian. To validate the randomized policy, we simulated the reference system by sampling, at each time-step, actions from the pf $g_{{u}}(u_k\mid\bv{x}_{k-1})$ and by then applying these actions to the dynamics in (\ref{eqn:pendulum_dynamics}).  In Figure \ref{fig:examples_simulation} the results of this process are illustrated.  The figure, which was obtained by performing $50$ simulations, clearly shows that the randomized policy is able to stabilize the pendulum around its unstable equilibrium point.

\begin{Remark}
See our github repository for all data, pfs and  code to replicate all the parts of the running example are given at the repository.
\end{Remark}

\section{Reinforcement Learning through the lenses of Problem \ref{prob:general_problem}}\label{sec:RL}
Throughout this section, we  do not assume that the agent knows the cost it is trying to minimize: it only receives a {\em cost signal} once an action is made \cite{RL_intro_2018,RL_a_survey,MLStory,Recht,offline_tuto,BUSONIU20188}. Typically, it is assumed that the cost signal received by the agent depends at each $k$ only on $\bv{X}_{k-1}$ and $\bv{U}_k$.  This leads to the following
\begin{Assumption}\label{asn:cost}
The cost in Problem \ref{prob:general_problem} is given by:
\begin{equation}\label{eqn:RL_generic_cost}
\begin{split}
   \mathbb{E}_{f}\left[c_{{1:T}}(\mathbf{X}_{0},\ldots,\bv{X}_T, \mathbf{U}_{1}, \ldots,\bv{U}_T)\right]  & = \mathbb{E}_{{f}}\left[\sum_{k\in1:T}c_{k}(\mathbf{X}_{k-1}, \mathbf{U}_{k})\right]\\
   & = -\mathbb{E}_{{f}}\left[\sum_{k\in1:T}d_kr_{k}(\mathbf{X}_{k-1}, \mathbf{U}_{k})\right],
   \end{split}
\end{equation}
where $d_k$ is a discount factor, $c_k(\cdot,\cdot)$ (resp. $r_k(\cdot,\cdot)$) is the cost (resp. reward) received at time-step $k$ by the agent when the state is $\bv{x}_{k-1}$ and $\bv{u}_k$ is applied. 
\end{Assumption}

{In (\ref{eqn:RL_generic_cost}) we used the shorthand notation $f$ to denote the pf $f(\boldsymbol{\Delta}_{0:T})$.} The sum in (\ref{eqn:RL_generic_cost}) is the cumulative reward obtained by the agent and Assumption \ref{asn:cost} is a standing assumption throughout the rest of the paper. In the RL terminology  $f_{\mathbf{u},k}(\mathbf{u}_{k}\mid \mathbf{x}_{k-1})$ is the {\em target policy}. The optimal target policy,  $f^\ast_{\mathbf{u},k}(\mathbf{u}_{k}\mid \mathbf{x}_{k-1})$, is the policy that the agent wishes to learn. Crucially, in certain RL algorithms the agent attempts to learn the target policy by following a different policy. It is this latter policy that is followed by the agent and determines its behavior.  For this reason, such a policy is termed as {\em behavior policy} and we denote it by  $\mu_{{k}}(\hat{\mathbf{u}}_{k}\mid \mathbf{x}_{k-1})$, where we are using the {\em hat} symbol to stress in our notation that the action generated by the behavior policy is different from the action that would have been obtained if the target policy was used.  
Target and behavior policies might depend on each other and this functional dependency can be expressed by adding a constraint to Problem \ref{prob:general_problem}. 
\begin{Remark}
Typical choices for the discount factor in (\ref{eqn:RL_generic_cost}) include: (i) constant, for example $d_k = 1/T$, $\forall k$; (ii) discounted, i.e., $d_k = \gamma^k, \gamma \in [0,1]$; (iii) myopic,  i.e., $d_0 = 1, d_k = 0, \forall k \ne 0$; (iv) final return, i.e.,  $d_T = 1, d_k = 0, \forall k\ne T$.
\end{Remark}
Based on these considerations, we formulate the following learning-oriented DM problem (simply termed as RL DM problem) derived from Problem \ref{prob:general_problem}:

\begin{Problem}\label{prob:RL_problem}
Find $\{f^{\ast}_{\bv{u},k}(\mathbf{u}_{k}\mid \mathbf{x}_{k-1})\}_{k \in 1:T}$  such that:

\begin{subequations}
    \begin{alignat}{3}
    \{f_{\bv{u},k}^{\ast}(\mathbf{u}_{k}\mid \mathbf{x}_{k-1})\}_{k \in 1:T} & \in \underset{\{f_{\bv{u},k}(\mathbf{u}_{k}\mid \mathbf{x}_{k-1})\}_{k \in 1:T}}{{\arg \min}}\  \mathbb{E}_{f}[c_{{1:T}}(\mathbf{X}_{0},\ldots,&&\bv{X}_T, \mathbf{U}_{1}, \ldots,\bv{U}_T)]\label{RLProb_cost}\\
    & s.t. \ \mathbf{x}_{k} \sim f_{\mathbf{x}{,k}}(\mathbf{x}_{k}\mid \mathbf{x}_{k-1}, \mathbf{u}_{k} = \hat{\mathbf{u}}_{k}),   && \ \forall k\in 1:T;\label{eqn:RLcst1}\\
    &\  \ \ \hat{\mathbf{u}}_{k} \sim \mu_{k}(\hat{\mathbf{u}}_{k}\mid \mathbf{x}_{k-1}), && \ \forall k\in 1:T;\label{eqn:RLcst2}\\
        &  \ \ \ \mathbf{u}_{k} \sim f_{\bv{u}{,k}}(\mathbf{u}_{k}\mid \mathbf{x}_{k-1}), && \ \forall k\in 1:T;\label{eqn:RLcst3}\\
    &\ \ \  \mathcal{G}\{f_{\mathbf{u}{,k}}(\mathbf{u}_{k}\mid \mathbf{x}_{k-1}),\mu_{k}(\hat{\mathbf{u}}_{k}\mid \mathbf{x}_{k-1})\}=0,  && \ \forall k\in 1:T;\label{eqn:RLcst4}\\
        &\ \ \  f_{\bv{u}{,k}}(\mathbf{u}_{k}\mid \mathbf{x}_{k-1}),\mu_{k}(\hat{\mathbf{u}}_{k}\mid \mathbf{x}_{k-1}) \in \mathcal{P}.\label{eqn:RLcst5}
    \end{alignat}
\end{subequations}
\end{Problem}
Problem \ref{prob:RL_problem} was obtained by relaxing the third, fourth and fifth constraint in Problem \ref{prob:general_problem}. The constraint set was also modified to take into account the presence of the behavior policy.  Constraint (\ref{eqn:RLcst1}) captures the fact that $\bv{x}_k$ is obtained by sampling from the probabilistic system description when the {\em previous} state is $\bv{x}_{k-1}$ and the action is obtained by sampling from the behavior policy rather than from the target policy (see also constraint (\ref{eqn:RLcst2})). The fact that the target and the behavior policies can be linked with each other {is formalized via the functional constraint (\ref{eqn:RLcst4})}: in the next sub-sections, we give a number of examples for this constraint and show how different choices lead to different exploration strategies. {Note  that, in some algorithms, there is no relationship between behavior  and  target policy. In these cases, constraint (\ref{eqn:RLcst4}) can be relaxed.  For example,  certain versions of  Q-Learning  make use of a behavior policy that corresponds to a random walk, see e.g., \cite{bomberman}. This can be embedded into Problem \ref{prob:RL_problem} by assuming that the behavior policy is the uniform pf.} 
\begin{Remark}\label{rem:formulation_RL}
The constraints of Problem \ref{prob:RL_problem} can also be formalized, see e.g., \cite{Recht,matni2019selftuning},  via stochastic difference equations (see the discussion in Section \ref{sec:pdfs} {and Remark \ref{rem:Wk}}). In this case, a common choice is to search for policies that are deterministic in the state (and in turn this leads to specify the control input via $\mathbf{U}_k = \pi_k(\mathbf{X}_{k-1})$). Rewards that contain some exogenous information
can be also considered.
In this case, in Problem \ref{prob:RL_problem} we have $r_k(\bv{X}_{k-1}, \bv{U}_k) := \E[R_k(\mathbf{X}_{k-1},\mathbf{U}_k, \mathbf{{W}}_k])$, where the expectation is taken over the pf from which $\mathbf{{W}}_k$ is sampled. { This notation highlights the presence of exogenous information  (see also the related discussion in Section \ref{sec:closed_loop_systems}).}
\end{Remark}
\begin{Remark}
When all the constraints except (\ref{eqn:RLcst1}) are relaxed and the dynamical systems formalism (see Remark \ref{rem:formulation_RL}) is used,  then Problem \ref{prob:RL_problem} becomes the one considered in \cite{Recht} to survey RL algorithms. In such a paper,  policy-based and value-based methods are surveyed through an  optimization framework.
\end{Remark}

\subsection{Tackling Problem \ref{prob:RL_problem}}\label{RLtaxo}

We now survey methods to solve Problem \ref{prob:RL_problem}. In doing so, we make use of the following:
\begin{Assumption}\label{asn:stationary}
The expectation in (\ref{eqn:RL_generic_cost}) is taken over the pf 
\begin{equation}\label{eqn:full_joint}
    f:= f(\mathbf{x}_0)\prod_{k\in1:T}f_\mathbf{u}(\mathbf{u}_{k}\mid \mathbf{x}_{k-1}) f_\mathbf{x}(\mathbf{x}_{k}\mid \mathbf{x}_{k-1}, \mathbf{u}_{k}).
\end{equation}
\end{Assumption}
Assumption \ref{asn:stationary} formalizes the fact that: (i) the probabilistic system description is stationary; (ii) the optimal solution is searched through target policies that are stationary.  
In turn, this implies that the decision variable in Problem \ref{prob:RL_problem} is $f_\mathbf{u}(\mathbf{u}_{k}\mid \mathbf{x}_{k-1})$. We organized the survey of the methods along three {\em dimensions}: (i)  model-based vs model-free; (ii) policy-based vs value-based; (iii) offline vs off-policy vs on-policy.  Essentially, the first dimension accounts for what the agent {\em knows} about the system, the second accounts for {\em how} the problem is solved and the third accounts for {\em when} data become available to the agent.  We now discuss each dimension and,  for each dimension, we give examples of state-of-the-art RL algorithms that fall in that dimension.

\subsubsection{Model-based vs model-free RL}

The first dimension arises from the knowledge of the probabilistic description of the system, i.e., from the knowledge of the pf in (\ref{eqn:RLcst1}).\\

\noindent {\bf Model-free reinforcement learning.} 
In a broad sense, {\em model-free} RL algorithms attempt to learn a policy without the knowledge of $f_\mathbf{x}(\bv{x}_k\mid\mathbf{x}_{k-1}, \mathbf{u}_{k})$.  Popular examples of model-free algorithms include Q-Learning \cite{qlearning} and SARSA \cite{SARSA}, which find a policy by building a state-action value function and picking the optimal action. Instead, REINFORCE \cite{reinforce} is a model-free RL algorithm which, rather than computing a value function,  learns a policy by estimating the gradients of the cost.  We refer readers to Section \ref{sec:policy_value} for a discussion on value-based and policy-based algorithms.  We also recall  certain model-free Actor-Critic algorithms \cite{AC1,AC2,ACsurvey} which perform learning by  building a policy and a set of value functions simultaneously (see also Section \ref{sec:policy_value}).\\

\noindent {\bf Model-based reinforcement learning.} When available,  $f_\mathbf{x}(\bv{x}_k\mid\mathbf{x}_{k-1}, \mathbf{u}_{k})$ can be  leveraged to improve the learning process. Algorithms that either make use of  $f_\mathbf{x}(\bv{x}_k\mid\mathbf{x}_{k-1}, \mathbf{u}_{k})$ or that learn this pf are termed as {\em model-based} RL algorithms \cite{morel,I2A,combo}. Some model-based algorithms also work on an estimate of the environment  via {\em world models}, or simulators \cite{DYNAQ,worldModel,mbve,MBpolicy}. These simulators leverage estimates of $f_\mathbf{x}(\bv{x}_k\mid\mathbf{x}_{k-1}, \mathbf{u}_{k})$ to e.g., generate imaginary data points.  Dyna-Q \cite{DYNAQ} is an example of model-based RL algorithm that can be thought of as a model-based variant of Q-Learning. Within Dyna-Q, interactions with the environment are augmented with additional data obtained via an estimate of $f_\mathbf{x}(\bv{x}_k\mid\mathbf{x}_{k-1}, \mathbf{u}_{k})$.  We note how the line of research in RL that makes use simulated data to learn a policy fits within the model-based classification. For example, in \cite{bansal2018chauffeurnet}, simulations are used to build models that are then leveraged to learn autonomous driving policies in real traffic conditions.  {Model-Based Value Expansion  \cite{mbve} and Model-Based Prior Model-Free  \cite{mbmf} also elaborate on this principle by iteratively learning a probabilistic dynamic model. {Certain models are also used within e.g.,  AlphaGo Zero \cite{alphago_zero,bertsekas2021lessons}. In this case, the agent is given a full environment model (specifically, the rules of the game of Go) prior to training and the model is in turn exploited for training a neural network. } Finally, we also report  the Probabilistic Inference for Learning COntrol  \cite{pilco} algorithm. This is a policy-based algorithm that explicitly accounts for model uncertainty to reduce learning biases related to flawed model estimates.}

\subsubsection{Policy-based vs value-based RL}
\label{sec:policy_value}
This dimension is related to how Problem \ref{prob:RL_problem} is solved.  By policy-based RL we refer to the set of techniques and algorithms that aim at directly finding a solution to Problem \ref{prob:RL_problem}, eventually assuming a parametrization of $f_\mathbf{u}(\bv{u}_k\mid \mathbf{x}_{k-1})$ and hence moving from an infinite-dimensional (functional) to a finite-dimensional optimization problem. Instead,  value-based RL finds the solution to Problem \ref{prob:RL_problem} indirectly, via a suitable {\em value} function.  We start to survey this latter approach. 

\begin{Remark}
The presence of the behavior policy does not play any role in value-based vs policy-based classification. Therefore, for notational convenience,  {constraint (\ref{eqn:RLcst2})} in Problem \ref{prob:RL_problem} is relaxed in this section. For the same reason, we omit specifying that the solution to Problem \ref{prob:RL_problem} needs to be a pf. 
\end{Remark}

\noindent {\bf Value-based reinforcement learning.} We start with introducing the following short-hand notation:
\begin{equation}\label{eqn:truncated_constraints}
\begin{split}
& \bv{C}_{a:b}(\bv{x},\bv{u}) := \\
& \left\{\bv{x}_{a-1} = \bv{x}, \ \bv{u}_a=\bv{u}, \ \mathbf{x}_{k} \sim f_{\mathbf{x},{k}}(\mathbf{x}_{k}\mid \mathbf{x}_{k-1}, \mathbf{u}_{k}), \mathbf{u}_k \sim f_{\mathbf{u},{k}}(\mathbf{u}_k\mid \mathbf{x}_{k-1}), \ k\in a:b \right\},
\end{split}
\end{equation}
to denote the set of constraints (\ref{eqn:RLcst1}), (\ref{eqn:RLcst3}) of Problem \ref{prob:RL_problem} between $k=a$ and $k=b$, with the additional condition that $\bv{x}_{a-1} = \bv{x}$ and $\bv{u}_a = \bv{u}$.  Value-based methods rely on finding the solution to Problem \ref{prob:RL_problem} via a value function. 
 
The so-called {\em state-action} value function \cite{qlearning,MLStory,RL_intro_2018} is defined as follows
\begin{align}
    &Q_{a\rightarrow b}(\mathbf{x},\mathbf{u}) := \underset{\{f_{\mathbf{u},{k}}(\mathbf{u}_k\mid \mathbf{x}_{k-1})\}_{a+1:b}}{\max} \mathbb{E}_{f_{a:b}}\left[\sum_{k\in a:b}d_{k}r_{k}(\mathbf{X}_{k-1}, \mathbf{U}_{k}) \mid \bv{C}_{a:b}(\bv{x},\bv{u}) \right],\label{eqn:q-fcn}
\end{align}
and determines the best  reward/cost value that the agent solving Problem \ref{prob:RL_problem} can achieve (in expectation) given $\bv{x}_{a-1} = \bv{x}$ and $\bv{u}_a = \bv{u}$.  In the above expression $f_{a:b}:= f(\bv{x}_{a-1},\bv{u}_a,\ldots, \bv{x}_{b-1},\bv{u}_b \mid \bv{x}_{a-2},\bv{u}_{a-1})$ is the pf of the evolution of the closed-loop system between $k=a$ and $k=b$. Following Bayes rule we have
\begin{equation*}
    \begin{split}
        f_{a:b} & = f(\bv{x}_{a-1},\bv{u}_a,\ldots, \bv{x}_{b-1},\bv{u}_b\mid \bv{x}_{a-2},\bv{u}_{a-1})\\
                & = f(\bv{x}_{a-1},\bv{u}_a\mid \bv{x}_{a-2},\bv{u}_{a-1})f(\bv{x}_{a},\bv{u}_{a+1},\ldots, \bv{x}_{b-1},\bv{u}_b\mid \bv{x}_{a-1},\bv{u}_a)\\
                & := f_a f_{a+1:b}
    \end{split}
\end{equation*}
This, together with the fact  that $\mathbf{u}$ and $\mathbf{a}$ are fixed,  is crucial to obtain an expression for $Q_{a\rightarrow b}(\mathbf{x},\mathbf{u})$ that enables its recursive computation. Indeed, from the definition of $Q_{a\rightarrow b}(\mathbf{x},\mathbf{u})$,  the following chain of identities (omitting explicitly writing the constraints) can be obtained ({without requiring Assumption \ref{asn:stationary}}):
\begin{equation}\label{eqn:Q_recursion}
\begin{split}
    Q_{a\rightarrow b}(\mathbf{x},\mathbf{u}) & =\E_{f_a}\left[d_{a}r_{a}(\mathbf{x},\mathbf{u}) +  \underset{\{f_{\mathbf{u},{k}}(\mathbf{u}_k\mid \mathbf{x}_{k-1})\}_{a+1:b}}{\max} \mathbb{E}_{f_{a+1:b}}\left[\sum_{k\in a+1:b}d_{k}r_{k}(\mathbf{X}_{k-1}, \mathbf{U}_{k})\right] \right]\\
    & = \E_{f_a}\left[d_{a}r_{a}(\mathbf{x},\mathbf{u}) +  \underset{\mathbf{u}'\sim f_{\bv{u},{a+1}}(\bv{u}_{a+1}\mid \bv{x}_a)}{\max} Q_{a+1\rightarrow b}(\mathbf{X}_{a},\mathbf{u}')\right]  \\
    & =\E_{f_a}\left[d_{a}r_{a}(\mathbf{x},\mathbf{u}) +  \underset{\mathbf{u}'}{\max} Q_{a+1\rightarrow b}(\mathbf{X}_{a},\mathbf{u}')\right].
\end{split}
\end{equation}
In the above expression, the first identity follows directly from the definition of $Q_{a\rightarrow b}(\mathbf{x},\mathbf{u})$, from the fact that from the definition of the constraints in (\ref{eqn:truncated_constraints}) we have $\bv{x}_{a-1} = \bv{x}$ and $\bv{u}_a=\bv{u}$ and from the fact that the decision variables do not depend on the pf over which the outer expectation is taken. The second identity follows again from the definition of value function and the last identity follows from the fact that $Q_{a+1\rightarrow b}(\mathbf{x}_{a},\mathbf{u}')$ depends directly on the control (and not the underlying pf). {The optimal value for the DM problem is given by $\underset{\bv{u}}{\max}  Q_{1 \rightarrow T}(\mathbf{x}_{0},\mathbf{u})$ which, following standard dynamic programming arguments, can be computed via backward recursion. From the same arguments, it also follows that, at each $k$:}
 \begin{align}\label{eqn:greedy_policy}
&f^\ast_{\mathbf{u}{,k}}(\mathbf{u}_{k}\mid \mathbf{x}_{k-1}) = \mathbbm{1}_{\bv{u}^\ast_{k}}(\bv{U}_k),\\
&{\bv{u}^\ast_{{k}} \in \underset{{\bv{u}}}{\arg\max}  Q_{{k}\rightarrow T}(\mathbf{x}_{{k-1}},\mathbf{u})}.
\end{align}
Computational barriers exist that prevent computing $Q_{1\rightarrow T}(\mathbf{x},\mathbf{u})$. In order to overcome these barriers,  different approximation techniques for $Q_{1\rightarrow T}(\mathbf{x},\mathbf{u})$ have been proposed under a wide range of technical conditions, see e.g.,  \cite{MLStory,qlearning,Powell_book,BUSONIU20188}.  {Perhaps, the most popular RL algorithm relying on these approximation methods is Q-Learning \cite{qlearning}, which will be further described in Section \ref{policy_line}.  A complete survey of} these approximation techniques goes beyond the scope of this paper and we refer to e.g., \cite{MC1} for a comprehensive monograph.
\begin{Remark}\label{rem:noisetypes}
The expectation in (\ref{eqn:q-fcn}) can be thought of as taken over all possible {\em types} uncertainties in $k\in a:b$ and this includes uncertainties on the state and control input.  Consider the case where: (i)  the state is generated by $\mathbf{X}_{k} = f(\mathbf{X}_{k-1}, \mathbf{U}_{k}, \mathbf{W}_k)$; (ii) one searches for policies that are deterministic in the state, i.e.,  $\mathbf{U}_k = \pi(\mathbf{X}_{k-1})$; (ii) the reward is given by $R_k(\bv{X}_{k-1},\bv{U}
_k,\bv{W}_k)$. Then the expectation in  (\ref{eqn:q-fcn}) needs only to be taken over the pfs from which $\{\bv{W}_k\}_{a:b}$ is sampled and from (\ref{eqn:Q_recursion}) the classic, see e.g., \cite[Chapter $11$]{MLStory}, recursion for the Q-function can be recovered.
\end{Remark}
The so-called {\em state} value function \cite{RL_intro_2018,MLStory} is defined  as follows:
\begin{align}\label{eqn:v_state_values}
    &V_{a\rightarrow b}(\mathbf{x}) = \underset{\{f_\mathbf{u}(\mathbf{u}_k\mid \mathbf{x}_{k-1})\}_{a+1:b}}{\max} \mathbb{E}_{f_{a:b}}\left[\sum_{k\in a:b}d_{k}r_{k}(\mathbf{x}_{k-1}, \mathbf{u}_{k})\mid  \bv{C}_{a:b}^v(\bv{x},\bv{u})\right],
\end{align}
where $\bv{C}_{a:b}^v$ is obtained by removing the constraint that $\bv{u}_a = \bv{u}$ from (\ref{eqn:truncated_constraints}). By definition, we  also have that
$V_{a\rightarrow b}(\mathbf{x}) = \underset{\bv{u}}{\max} Q_{a\rightarrow b}(\mathbf{x},\mathbf{u})$.
Recursive equations analogous to (\ref{eqn:Q_recursion}) can be obtained for the state value function, see e.g., \cite[Chapter $3$]{RL_intro_2018}. These lead, in particular, to temporal difference algorithms \cite{Sutton:1988up,5152964}. Finally,  certain RL algorithms combine the state-action and the state functions by defining the {\em advantage} function \cite{10.5555/3045390.3045601,baird_aa}: 
$A_{a\rightarrow b}(\bv{x},\bv{u}) := Q_{a\rightarrow b}(\mathbf{x},\mathbf{u}) - V_{a\rightarrow b}(\mathbf{x})$.
The advantage function can be estimated via two separate estimators, i.e., one for each of the value functions and the main benefit of using this architecture is the possibility of generalizing learning across actions without imposing any change to the underlying RL algorithm. {Advantage Actor-Critic algorithms, or A2C \cite{AC2}, leverage this advantage function to improve learning.}
\begin{Remark}
{The functions discussed above can either be represented by a table,  or parametrized and directly approximated.  In the latter case, when the parameters are the weights of a deep neural network, the methods fall under the label of {\em deep} RL algorithms, see e.g., \cite{A3c, SAC, DeepDynaQ}. Popular tabular methods are  SARSA and Q-Learning, which are described  in Section \ref{policy_line}. }
\end{Remark}

\noindent {\bf Policy-based reinforcement learning.} In policy-based reinforcement learning, the target policy is found without passing through the computation of the value function.  Within these methods, the policy often has a fixed structure. That is, in the context of Problem \ref{prob:RL_problem}, the optimization is performed over a family of pfs, parametrized in a given vector of parameters.  In turn, this is equivalent to restrict the feasibility domain of Problem \ref{prob:RL_problem} by changing its last constraint (\ref{eqn:RLcst5}) so that  $f_{\bv{u}}(\bv{u}_k \mid \mathbf{x}_{k-1})$ does not just belong to $\mathcal{P}$ but rather to a parametrized pf family. A Typical choice for the parametrized family are exponential families; another approach is that of using neural networks to parametrize the pf. In this case,  the vector parametrizing the pf are  the weights of the network.  In what follows, we denote the vector parametrizing the pf by $\boldsymbol{\theta}$ and, to stress the fact that the policy is parametrized, we write $f_\mathbf{u}(\bv{u}_k \mid \mathbf{x}_{k-1}) = f_{\mathbf{u},\boldsymbol{\theta}}(\bv{u}_k \mid \mathbf{x}_{k-1})$. Hence,  the decision variable in Problem \ref{prob:RL_problem} becomes $\boldsymbol{\theta}$ and the goal is that of finding the optimal $\boldsymbol{\theta}^\ast$.  With this formulation, a possible approach to solve Problem \ref{prob:RL_problem} relies on estimating the gradient of the objective function, see e.g., \cite{MLStory,offline_tuto}.  As discussed in \cite{offline_tuto}, where $d_k = d^k$,  a rather direct estimate relies on expressing the gradient of the objective w.r.t.  $\boldsymbol{\theta}$ as follows:
\begin{equation}\label{eqn:gradient_base}
    \nabla_\theta \left(\mathbb{E}_{{f}}\left[\sum_{k\in1:T}d^kr_{k}(\mathbf{X}_{k-1}, \mathbf{U}_{k})\right]\right) = \mathbb{E}_{f}\left[\sum_{k\in1:T}d^k\nabla_\theta \ln f_{\mathbf{u},\theta}(\mathbf{u}_{k}\mid \mathbf{x}_{k-1})\hat{A}(\bv{x}_{k-1},\bv{u}_k)\right],
\end{equation}
where $\hat{A}$ is a return estimate, which can be obtained via e.g., Monte-Carlo methods.  Algorithms such as REINFORCE {\cite[Chapter $12$]{MLStory}}, \cite{reinforce} attempt to find $\boldsymbol{\theta}^\ast$ by sampling from $ f_{\mathbf{u},\theta}(\bv{u}_k \mid \mathbf{x}_{k-1})$ to build an estimate of the gradient and hence running gradient ascent iterates. An alternative to use Monte-Carlo methods is to estimate $\hat{A}$ via a separate neural network (i.e., the critic) which updates the policy alongside with the actor. As the critic network is essentially a value estimator, these {actor-critic} methods essentially combine together value-based and policy-based iterates \cite{AC1,AC2}. Actor-critic methods can achieve better convergence performance than pure critic algorithms \cite{AC1} and a further evolution of these algorithms relies on parallel training \cite{A3c} where several actors and critics are trained simultaneously.  In this context, we also recall  the soft actor critic algorithm \cite{SAC}, which encourages randomness in policies by regularizing the objective function with an entropy term.

\subsubsection{On-policy,  off-policy and offline RL}\label{policy_line}

This dimension accounts for when the data used to find the policy solving Problem \ref{prob:RL_problem} are collected.  The agent can indeed use either the target policy, i.e.,  $f_{\bv{u}}(\bv{u}_k\mid \mathbf{x}_{k-1})$ in Problem \ref{prob:RL_problem}, to collect the data or a suitably defined behavior policy, i.e., $\mu(\hat{\bv{u}}_k\mid \mathbf{x}_{k-1})$ in Problem \ref{prob:RL_problem}, to encourage exploration. This data collection process can be online or, as an alternative, data can be all collected before runtime via some potentially unknown behavior policy. \\

\noindent {\bf On-policy reinforcement learning.} On-policy (also known as fully online) RL algorithms rely on collecting data online via the target policy. This means that, in Problem \ref{prob:RL_problem},  constraint (\ref{eqn:RLcst4}) becomes $f_{\mathbf{u}}(\mathbf{u}_{k}\mid \mathbf{x}_{k-1})=\mu(\mathbf{u}_{k}\mid \mathbf{x}_{k-1})$.  A classic example of on-policy RL algorithm is SARSA \cite{SARSA}, which is a {tabular} value-based algorithm. SARSA aims to estimate the state-action function by using a target policy that is derived, at each $k$, by the current estimate of the Q-table (i.e., the tabular representation of the Q-function). {In the discounted case, with {infinite time horizon and} stationary rewards,} this estimate is updated at each $k$ as follows:
\begin{align}\label{eqn:SARSA}
    &{Q}_{new}(\mathbf{x}_{k-1}, \mathbf{u}_k) \gets {Q}_{old}(\mathbf{x}_{k-1}, \mathbf{u}_k) + \alpha (R_k + \gamma \mathcal{Q}_{old}(\mathbf{x}_{k}, \mathbf{u}_{k+1}) - \mathcal{Q}_{old}(\mathbf{x}_{k-1}, \mathbf{u}_{k})),
\end{align}
where $\alpha \in (0,1)$ is a learning rate,  $R_k$ is the reward signal {received} by the agent when $\bv{u}_k$ is selected and the system is in state $\bv{x}_{k-1}$.  In SARSA, the element $(\bv{x}_{k-1},\bv{u}_k)$ is updated based on an action (i.e., the target action) that is obtained from the Q-table.  SARSA (as well as Q-Learning, discussed within the off-policy methods) can be improved by sampling several actions and states from the Q-table.  This is done to improve   the estimates of the Q-table and, in the extreme case where a whole episode is played before updating the policy, this technique becomes a Monte-Carlo algorithm \cite{MC1a,MC2a}. The choice of how $\bv{u}_k$ is selected from the Q-table has a key impact on the algorithm performance. Using a greedy policy (\ref{eqn:greedy_policy}) decreases exploration and can potentially prevent from learning optimal actions \cite{RL_intro_2018}. To mitigate this,  randomness can be added to the greedy policy, thus obtaining a $\varepsilon$-greedy policy. Another popular choice is the so-called softmax policy \cite{mellowmax},  which, by defining  $f_{\bv{u}}(\bv{u}_k \mid \mathbf{x}_{k-1}) $ as a Boltzmann pf \cite{boltzmann}, adds a design parameter to the algorithm. This is the {temperature}: if the temperature is $0$, then the agent behaves greedily, while increasing it encourages exploration. While SARSA is perhaps the most popular on-policy RL algorithm, we recall  other algorithms such as {the gradient-based Trust Region Policy Optimization  \cite{trpo}  algorithm, which uses a KL-divergence (see Definition \ref{def:KLdiv}) constraint between the starting and updated policy at each iteration,  and its approximate counterpart, Proximal Policy Optimization \cite{ppo}.}\\


\noindent {\bf Off-policy reinforcement learning.} In {\em off-policy} RL the behavior policy is different from the target policy. That is, the policy that the agent tries to learn is different from the policy actually used by the agent.  As a result,  the target policy is learned from trajectories sampled from the behavior policy. The functional relationship between target and behavior policy is expressed via constraint (\ref{eqn:RLcst4}) of Problem \ref{prob:RL_problem}.  Q-Learning is perhaps the most popular off-policy RL algorithm \cite{qlearning,human_lvl_control,doubleQL}.  As SARSA, this is a tabular value-based RL algorithm, which is based on the use of the state-action function (i.e., the Q-table).  However,  differently from SARSA,  Q-Learning updates its estimate of the state-action function via a greedy policy, while the agent behavior is determined by using the behavior policy (that encourages exploration). The resulting update rule is {(again in the discounted case {with stationary reward and infinite-time horizon})}:
\begin{align}\label{eqn:Qlearning}
    &{Q}_{new}(\mathbf{x}_{k-1,} \mathbf{u}_k) = {Q}_{old}(\mathbf{x}_{k-1}, \mathbf{u}_k) + \alpha (R_k + \gamma \underset{\mathbf{u}}{\max}{Q}_{old}(\mathbf{x}_{k}, \mathbf{u}) - {Q}_{old}(\mathbf{x}_{k-1}, \mathbf{u}_{k})),
\end{align}
where $\alpha$ and $R_k$ are defined as in (\ref{eqn:SARSA}).  The key difference between  the update rule in (\ref{eqn:Qlearning}) and the one in (\ref{eqn:SARSA}) is that in the former case ${Q}_{new}(\mathbf{x}_{k-1,} \mathbf{u}_k) $ depends on $\underset{\mathbf{u}}{\max}{Q}_{old}(\mathbf{x}_{k}, \mathbf{u})$.
It is interesting to note how the behavior policy is defined.  A typical choice in Q-Learning is to pick  the behavior policy as a $\varepsilon$-greedy version of the target policy. This choice for the behavior policy formalizes the fact that agent randomly explores non-greedy actions with some design probability $\varepsilon$. {In turn, the link between the behavior and the target policy can be captured} via constraint (\ref{eqn:RLcst4}) of Problem \ref{prob:RL_problem}. { Indeed, $\mu(\hat{\mathbf{u}}_k\mid \mathbf{x}_{k-1})$ can be written as:
\begin{equation}
\mu(\hat{\mathbf{u}}_k\mid \mathbf{x}_{k-1}) = (1-\epsilon)\cdot f_{\mathbf{u}}(\mathbf{u}_{k}\mid \mathbf{x}_{k-1}) + \epsilon \cdot \text{unif}(\bv{U}_k).
\end{equation}
In the above expression, which can be formalized via (\ref{eqn:RLcst4}), $f_{\mathbf{u}}(\mathbf{u}_{k}\mid \mathbf{x}_{k-1})  = \mathbbm{1}_{\tilde{\bv{u}}_{k}}(\bv{U}_k)$, with $\tilde{\bv{u}}_{{k}} \in \underset{{\bv{u}}}{\arg\max}  Q_{old}(\mathbf{x}_{{k-1}},\mathbf{u})$, and $\text{unif}(\bv{U}_k)$ denoting the uniform pf over the action space.}
 The greediness does not necessarily need to be constant over time and a common choice is indeed that of decreasing $\varepsilon$ gradually over the episodes ({this is for example done in} \cite{human_lvl_control}).  Another choice to relate the behavior policy and the target policy include the use of the softmax  policy. {This policy is  given by the Boltzmann pf ${\mu}({\hat{\mathbf{u}}}_k\mid \mathbf{x}_{k-1}) = \frac{e^{{Q_{{old}}}(\mathbf{x}_{k-1},\mathbf{u}_k)/\rho}}{\sum_{\mathbf{x}\in\mathcal{X}}e^{{Q_{{old}}}(\mathbf{x},\mathbf{u}_k)/\rho}}$ where $\rho$ is the {temperature}.} {The link between the behavior policy and the target policy can be again captured via  (\ref{eqn:RLcst4}). Indeed, the behavior policy can be written as:
\begin{equation}
    \mu(\hat{\mathbf{u}}_k\mid \mathbf{x}_{k-1}) = \frac{f_{\mathbf{u}}(\mathbf{u}_{k}\mid \mathbf{x}_{k-1})\cdot e^{{Q_{old}}(\mathbf{x}_{k-1},\mathbf{u}_k)/\rho} + (1-f_{\mathbf{u}}(\mathbf{u}_{k}\mid \mathbf{x}_{k-1}))e^{{Q_{old}}(\mathbf{x}_{k-1},\mathbf{u}_k)/\rho}}{\sum_{\mathbf{x}\in\mathcal{X}}e^{{Q_{old}}(\mathbf{x},\mathbf{u}_k)/\rho}},
\end{equation}
with $f_{\mathbf{u}}(\mathbf{u}_{k}\mid \mathbf{x}_{k-1})  = \mathbbm{1}_{\tilde{\bv{u}}_{k}}(\bv{U}_k)$, $\tilde{\bv{u}}_{{k}} \in \underset{{\bv{u}}}{\arg\max}  Q_{old}(\mathbf{x}_{{k-1}},\mathbf{u})$.}
Q-Learning can also be implemented via function approximators \cite{human_lvl_control,mnih2013playing}, i.e.,  neural networks to approximate Q-functions. These algorithms go under the label of {\em deep} RL: popular algorithms are C51 \cite{C51} and {QR-DQN} (Quantile Regression DQN) \cite{qrdqn}. \\

\noindent {\bf Offline reinforcement learning.} In offline RL, see \cite{offline_tuto} for a detailed survey, data are collected once, before runtime,  via an offline behavior policy \cite{pseudoDynaQ}.  The behavior policy,  $\mu(\hat{\bv{u}}_k \mid \mathbf{x}_{k-1})$ in Problem \ref{prob:RL_problem} can be unknown to the agent.  As noted in \cite{pessimism},  where a link between offline RL and imitation learning is unveiled, two types of offline {\em data collection} methods have shown promising empirical and theoretical success: data collected through expert actions and uniform coverage data. In the former case,   the offline behavior policy can be thought of  as an {\em expert} that illustrates some desired behavior to fulfill the agent task. Then, offline RL intersects with the {\em imitation learning} framework \cite{imitationLearning, 9317713, pessimism}. In the latter case,  instead, the goal of the offline behavior policy is that of widely covering the state and action spaces to get informative insights \cite{D4RL,WILDS}.  In both cases, the functional relation between the target and the behavior policy is still captured via constraint (\ref{eqn:RLcst4}) of Problem \ref{prob:RL_problem} and the same considerations highlighted above for off-policy RL still hold in this case.  However, the key conceptual difference is that now the behavior policy has the role of {\em grounding} the target policy. In this context, a key challenge is encountered when the agent meets out-of-distribution samples (see e.g., \cite{offline_tuto,mopo}). It is indeed known that a discrepancy between the distributions over states and actions sampled from the behavior and target policy can negatively impact performance see e.g., \cite{CQL,noexplo,lee2021representation} and references therein.  Interestingly,  \cite{wang2021what} has shown that sample efficient offline RL is not possible unless there is a low distribution shift (i.e., the offline data distribution is close to the distribution of the target policy), while \cite{lee2021addressing}  proposed to mitigate this issue by keeping separate offline and online replay buffers. We refer readers to these papers and to \cite{offline_tuto,pessimism} for a detailed survey on offline RL methods and mitigation strategies for the distribution shift.

\subsection{Comments}\label{RLcomments}
We give here a number of comments, which are transversal to the dimensions presented above.\\

\noindent {\bf Hie{r}archical algorithms.} Feudal, or hierarchical RL \cite{feudal1993}, is based on the idea of splitting the agent's algorithm into a high level decision-making and a low level control component. The lower level control guarantees the fulfillment of the goal set by the decision-making algorithm.  In particular, the so-called {\em manager} is trained to find advantageous directions in the state space, while the {\em worker} is trained to output the correct actions to achieve these directions.  In \cite{deepFeudal}, the manager is trained with a modified policy gradient algorithm (where the modification exploits the knowledge that the output is a direction to be fed to a worker), while the worker uses an advantage actor-critic algorithm. Other works include \cite{options_feudal,states_as_goals},  where the manager outputs sub-policies with termination conditions, and is only queried again when such conditions are met.\\

\noindent {\bf Safety.}  We refer to \cite{brunke2021safe} for a detailed review of learning-based control. Essentially, safety constraints can be embedded in Problem \ref{prob:RL_problem} in three ways, which in turn correspond to different safety levels \cite{brunke2021safe}. Namely, safety can be embedded by: (i)
 {\em encouraging} safety constraints -- this can be typically done by adding a penalty to the cost for action-state pairs that are considered {\em unsafe} by the designer. This can be handled in Problem \ref{prob:RL_problem} by including to the cost and an additional penalty; (ii) guaranteeing safety constraints {\em in probability} -- this can be achieved by adding box (or chance constraints) to Problem \ref{prob:RL_problem}. These constraints can be handled by bringing back in the formulation of Problem \ref{prob:RL_problem} constraints (\ref{genProb_ineqConstr}) and (\ref{genProb_safetySet}) of Problem \ref{prob:DDC_problem}; (iii) foreseeing {\em hard} constraints, adding them to the formulation of Problem \ref{prob:RL_problem}. Formally, these constraints can be handled by adding to the formulation constraint (\ref{genProb_safetySet}) of Problem \ref{prob:DDC_problem} with $\varepsilon_k = 0$. \\

\noindent {\bf Distributional and memory-based RL.} Distributional RL leverages ideas from  distributionally robust optim{iz}ation in order to maximize worst case rewards and/or to protect the agent against environmental uncertainty.  This leads to a  functional formulation of Bellman's optimality, where maximization of a stochastic return is performed with the constraint that the policies belong to a given subset (in the space of probability densities).  We refer readers to e.g.,  \cite{C51,distrib1,distrib2} which leverage the distributionally robust optimization framework in a number of learning applications.  We also mention  memory-based RL, which improves data efficiency by {\em remembering} profitable past experiences.  In e.g., \cite{memory_RL} the agent uses a memory buffer through which past, surprisingly profitable (as in, collecting more reward than expected a priori) experiences are stored. This experience is then leveraged when the agents encounters one of the states in this memory buffer: when this happens, the agent behavior policy is biased towards reiterating past decisions.  See e.g., \cite{memory_survey} for a survey of the memory-based RL.\\

{\noindent {\bf Links between RL and inference.} An interesting connection exists between RL and inference of probabilistic models. As noted in e.g., \cite{inference_joe,inference_levine}, by recasting the problem of computing decisions as an inference problem, one (besides unveiling an intriguing  duality between DM and inference) gets access to a number of widely established inference tools that can be useful when designing the algorithm. The broad idea is to introduce a binary random variable sometimes termed as {\em belief-in-optimality} \cite{inference_joe}, say $\bv{Z}_k$, with $\bv{z}_k = 1$ indicating optimality. The pf $p(\bv{z}_k=1\mid\bv{x}_{k-1},\bv{u}_{k})$ is usually chosen so that $p(\bv{z}_k=1\mid\bv{x}_{k-1},\bv{u}_{k})\propto \exp(\rho r_k(\bv{x}_{k-1},\bv{u}_k))$, where $\rho$ is the temperature parameter. The temperature parameter can be fixed using heuristics, such as being a multiple of the reward standard variation, when available, \cite{inference_policySearch}. When the model is known, the parameters of this distribution can be computed through a backward recursion, while when the plant is unknown they can be estimated through an update rule similar in spirit to Q-Learning. With this set-up, stationary randomized control policies can be computed by inferring the {\em messages} $p(\bv{z}_{k:T}\mid\bv{x}_{k-1},\bv{u}_k)$ and $p(\bv{z}_{k:T}\mid\bv{x}_{k-1})$, corresponding to the probability of a path being optimal when starting from a state (and, for the first message, by taking an action). Moreover, as noted in \cite{inference_levine},  the problem of inferring these messages can be recast as a divergence minimization problem. In turn, this problem is equivalent to maximizing the sum of rewards plus an entropy term. Also, in \cite{inference_rl_control}, a reformulation of a classic stochastic control problem in terms of divergence minimization is presented. Then, a direct application of Bayesian inference is used to iteratively solve the problem and model-free settings are considered. We note that additional specifications on the DM  problem can lead to complementary ways of leveraging the $\bv{Z}_k$'s to find optimal policies and this is leveraged for example in \cite{inference_LQG} for  Linear Quadratic Gaussian problems. Finally, we also recall that approximate Bayesian inference has also been used for advancing trajectory optimization in the context of covariance control \cite{https://doi.org/10.48550/arxiv.2103.06319} 
}

\begin{figure}[thbp]
\centering
\includegraphics[width=0.48\linewidth]{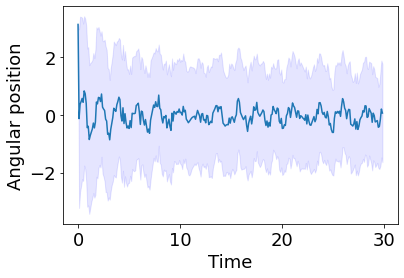}
\includegraphics[width=0.48\linewidth]{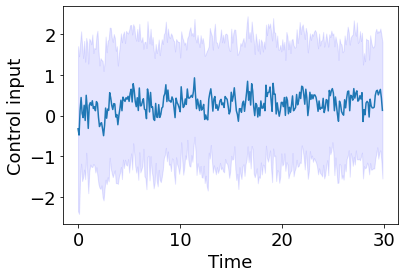}\\
\includegraphics[width=0.48\linewidth]{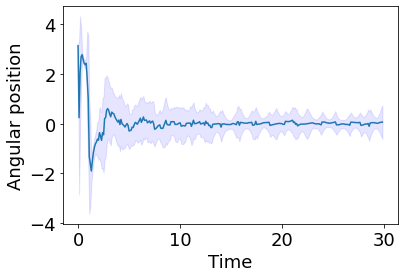}
\includegraphics[width=0.48\linewidth]{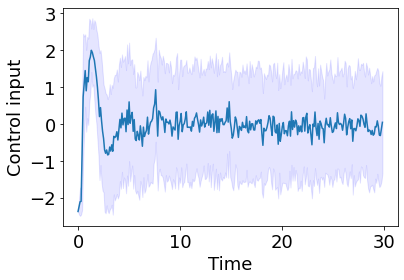}
\caption{swing-up via Q-Learning at the $2000$ (top panels) and $100000$ (bottom panels) checkpoints. Bold lines denote {means}, shaded areas the confidence intervals corresponding to the standard deviation. Figure obtained from $50$ simulations.}
\label{fig:examples_simulation_QL_2}
\end{figure}

\begin{myexp*}\normalfont We now stabilize the upward, unstable, equilibrium position of the pendulum  we previously introduced via the popular Q-Learning algorithm ({we refer interested readers to our github for the implementation details and for the stand-alone code}). We recall that Q-Learning is an off-policy, model-free, value-based algorithm that leverages a tabular representation of the state-action value function. Following the tabular nature of the algorithm, the first step to apply Q-Learning for the control of the pendulum (in the experiments, we considered a mass of {$1$}kg) was that of discretizing the action/state spaces. We used the same discretization that we presented when we obtained the probabilistic descriptions, i.e, $\mathcal{X}$ was discretized in a grid of $50\times50$ uniform bins and $\mathcal{U}$ was discretized in $20$ uniform bins. The reward signal received by the agent at each $k$  was instead given by $R_k = - \theta_k^2 - 0.1\omega_k^2$. The algorithm parameters were as follows: the learning rate was $\alpha = 0.5$, the discount factor was $\gamma = 0.99$ and the behavior policy was $\epsilon$-greedy with $\epsilon = 0.9$. Training episodes were $500$ time-steps long and the pendulum was set to the downward position at the beginning of each episode.

Given this set-up, we trained the Q-Learning agent using the behavior policy and a backup of the Q-table was performed at the following {\em checkpoints}: $20$, $200$, $2000$, $20000$ and $100000$ training episodes. This was done in order to obtain a {\em snapshot} of the policy learned by the agent as the number of training episodes increases. At each checkpoint, we controlled the pendulum using the policy learned by the agent. The policy was evaluated by running $50$ evaluation episodes, with each episode now being $300$ time-steps long and using the trained, greedy, policy. The result of this process is illustrated in Figure \ref{fig:examples_simulation_QL_1}, where the mean reward, together with its confidence interval, is reported. The mean reward was obtained by averaging the rewards obtained by the agent across all the evaluation episodes in the last  $100$ time-steps. Interestingly, in the figure it is shown that the agent learns to perform the task by first improving the mean reward then its confidence interval.

\begin{figure}[thbp]
\centering
\includegraphics[width=0.9\linewidth]{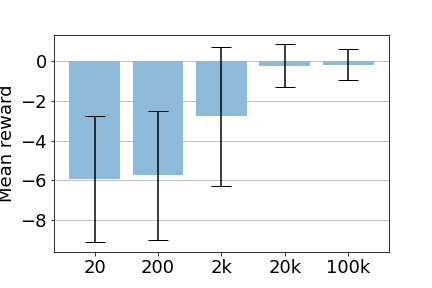}
\caption{mean reward over the last $100$ simulation steps of the evaluation episodes at each of the checkpoints of $20$, $200$, $2000$, $20000$ and $100000$ training episodes. The height of the  bars denotes the value of the mean reward achieved by the agent and the black lines represent the confidence intervals corresponding to the standard deviation.}
\label{fig:examples_simulation_QL_1}
\end{figure}

In order to further illustrate the behavior of the agent as the number of training episodes increases, we stored the evaluation results at each checkpoint. In Figure \ref{fig:examples_simulation_QL_2} the behavior of the pendulum is shown when the policies learned at $2000$ and $100000$ episodes are used to control the pendulum. In the figure it is clearly shown that, as the number of episodes increases, the agent learns to swing-up the pendulum (note the different confidence intervals for $\theta_k$ at $2000$ and $100000$  training episodes and how these are in agreement with Figure \ref{fig:examples_simulation_QL_1}).
\end{myexp*}

{
\subsection{Multi-armed bandits}\label{sec:MAB}
The framework considered in this paper also applies to {multi-armed bandit} problems. While we refer readers to e.g., \cite{lattimore2020bandit,slivkins2019introduction} for detailed monographs on this topic, we give here a discussion on this class of sequential DM problems through the lenses of Problem \ref{prob:RL_problem}. We recall that, throughout Section \ref{sec:RL}, we do not assume that the agent knows the cost and, for our discussion, we formulate the following DM problem derived  from Problem \ref{prob:RL_problem}:
\begin{Problem}\label{prob:multiArmed_problem}
Let $f:=f(\boldsymbol{\Delta}_{0:T})$. Find $\{f^{\ast}_{\bv{u},k}(\mathbf{u}_{k}\mid \mathbf{x}_{k-1})\}_{k \in 1:T}$  such that:
\begin{subequations}
    \begin{alignat}{3}
    \{f_{\bv{u},k}^{\ast}(\mathbf{u}_{k}\mid \mathbf{x}_{k-1})\}_{k \in 1:T} & \in \underset{\{f_{\bv{u},k}(\mathbf{u}_{k}\mid \mathbf{x}_{k-1})\}_{k \in 1:T}}{{\arg \max}}\  \mathbb{E}_{f}&&\left[\sum_{k \in 1:T}\mathbb{E}_{\mathbf{W}_k}[R_k(\mathbf{X}_{k-1}, \mathbf{U}_k,\mathbf{W}_k)]\right]\label{multiArmed_prob_cost}\\
  s.t. \ &\mathbf{x}_{k} \sim f_{\mathbf{x},k}(\mathbf{x}_{k}\mid \mathbf{x}_{k-1}, \mathbf{u}_{k}),   && \ \forall k\in 1:T;\label{eqn:MAcst1}\\
        &   \mathbf{u}_{k} \sim f_{\bv{u},k}(\mathbf{u}_{k}\mid \mathbf{x}_{k-1}), && \ \forall k\in 1:T;\label{eqn:MAcst2}\\
        & f_{\bv{u},k}(\mathbf{u}_{k}\mid \mathbf{x}_{k-1})\in \mathcal{P}.\label{eqn:MAcst3}
    \end{alignat}
\end{subequations}
\end{Problem}
Problem \ref{prob:multiArmed_problem} was obtained by reformulating Problem \ref{prob:RL_problem} as a reward maximization problem and by relaxing the constraints related to the behavior policy. The notation in (\ref{multiArmed_prob_cost}) is used to highlight the presence of exogenous information and the expectation inside the sum is taken, at each $k$, over the pf from which $\bv{W}_k$ is sampled. Also, the pfs $f_{\bv{x},k}(\mathbf{x}_{k}\mid \mathbf{x}_{k-1}, \mathbf{u}_{k})$ are not necessarily known by the agent. In this context, we note that a crucial aspect when modeling multi-armed bandit problems is the definition of state. The state variable can indeed be used to introduce the \emph{context} into the bandit problem. In this case, the DM problem is termed as contextual bandit problem and the pf from which the context is sampled is typically not known by the agent \cite[Chapter $12$]{MLStory}. Contextual bandits are a convenient framework to e.g., design context-aware recommender systems \cite{10.1007/978-3-642-34487-9_40} and embed users' feedback in an optimization process \cite{osti_1798214}; we refer to e.g., \cite{NIPS2007_4b04a686,pmlr-v15-chu11a,pmlr-v32-agarwalb14} for a number of technical results (and algorithms) on contextual bandits. We also note that a part of literature on multi-armed bandits includes in the definition of state a notion of \emph{belief}, which typically is an estimate of the reward distributions maintained by the agent \cite{powell2020state} (as such, the dynamics governing the evolution of the belief can be known to the agent). {Belief states from the literature include \cite[Chapter $1$]{Powell_book} the empirical mean reward over past episodes and parameters characterizing the posterior distribution of the rewards. These definitions of belief are also used within the resolution methods surveyed next.}
In a broader context, as noted in  \cite{powell2020state}, resolution approaches for Problem \ref{prob:multiArmed_problem} that leverage the use of {\em Gittins indices} are based on solving Bellman's equation where the state variable is the belief \cite{alma9960082294202441}. In Problem \ref{prob:multiArmed_problem}, the decision variables are randomized policies (the last constraint in the problem is a normalization constraint guaranteeing that the solution is a pf). We recall that, while policies that are deterministic in the state can be formally expressed as pfs, Problem \ref{prob:multiArmed_problem} also allows to consider situations, naturally arising in the Bayesian framework of e.g., \cite{JMLR:v13:may12a}, in which the optimal policy is randomized. 
Finally, we note that Problem \ref{prob:multiArmed_problem} captures bandits with both continuous and discrete state/action spaces. While a large body of the literature studies multi-armed problems with discrete action spaces \cite[Part VII]{lattimore2020bandit}, we refer readers interested in multi-armed problems with continuous actions to e.g., \cite{continuous_actions}. Given this set-up, we now discuss a number of resolution algorithms.\\
\newline
\noindent \textbf{Tackling Problem \ref{prob:multiArmed_problem}.} A key feature of the resolution methods, which aim at finding policies with guaranteed regret bounds \cite{lattimore2020bandit,slivkins2019introduction}, is that the agent attempts to solve Problem \ref{prob:multiArmed_problem} by building (and iteratively improving) an estimate of the reward  associated to each action. Certain algorithms, such as the  explore-then-commit algorithm, do this by explicitly enforcing an exploration phase (see below and e.g., \cite{ETC_review}). Instead, other algorithms directly embed a mechanism that favors exploration of less-used actions. This is the case of the Upper Confidence Bound algorithm (UCB); see e.g., \cite{UCB} and \cite[Chapter $8$]{slivkins2019introduction}. Specifically, for each pair of action and context (if any) UCB algorithms form an upper bound with a given confidence $\delta$, say $\bar{B}_k(\bv{x}_{k-1},\bv{u}_k)$, of the mean reward computed from the empirical mean reward (i.e., the belief state). Given this bound, the resulting policy is then of the form
\begin{equation}
    \begin{aligned}
    f_{\bv{u},k}(\mathbf{u}_{k}\mid \mathbf{x}_{k-1}) &= \mathbbm{1}_{\mathbf{u}_k}(\mathbf{U}_k),\\
    \mathbf{u}_k &= \underset{\mathbf{u}}{\arg\max}\ \bar{B}_{k}(\bv{x}_{k-1},\bv{u}).
    \end{aligned}
\end{equation}
One can obtain expressions for $\bar{B}_k(\bv{x}_{k-1},\bv{u}_k)$ that are inversely proportional to the number of times a given action was tried for a given context before time-step $k$ (for non-contextual bandits, the bound depends on the number of times each action is taken). Hence, with these bounds, exploration of less-used actions is favored. See also \cite{Auer_2010,10.5555/944919.944941} for more details on UCBs and the related principle of optimism in the face of uncertainty. The explore-then-commit algorithm  foresees an exploration-only phase with a subsequent commit phase within which the agent follows a greedy policy. The greedy policy selects the action with the highest empirical mean reward and in \cite{MLStory} it is noted how purely greedy algorithms without an exploration phase can  perform well for specific classes of problems (e.g., for linear bandits the bound on the regret is proportional to the square root of the time horizon \cite{greedy_bandits}) although in the worst case can incur in linear regret. Exploration can be either performed by selecting random actions \cite[Chapter $12$]{MLStory}, \cite[Chapter $8$]{slivkins2019introduction} or by choosing each action a given (pre-defined) amount of times \cite[Chapter $6$]{lattimore2020bandit}. We also recall the successive elimination algorithm \cite{slivkins2019introduction, MLStory}, which maintains an upper and lower bound on the average reward (again computed from the empirical mean reward) given by each action. Actions are \emph{deactivated} when their upper bound becomes lower than the lower bound of any other action. In between eliminations, exploration is carried out only on the actions that are still  active. \\
The above algorithms, originally introduced for non-contextual bandits, perform well also for contextual bandits with small context spaces but suffer a drop in performance as the context space becomes larger \cite[Chapter $8$]{slivkins2019introduction}, \cite[Chapter $18$]{lattimore2020bandit}. With respect to this aspect, we recall the linear UCB algorithm (LinUCB) which tackles large context spaces for linear bandit problems, i.e., problems for which the reward is a linear combination of \emph{features}. The algorithm  maintains a confidence region for the coefficients of the linear combination and, based on this region, calculates UCBs for the expected reward. Another idea, reminiscent of the explore-then-commit algorithm, is that of using (after the exploration phase) a classifier to identify the best policy \emph{in hindsight} (i.e., the policy that would have yielded the best reward during the exploration phase) by assigning to each context the action that obtained the best reward (computed from the belief state). This policy is then used by the agent. See \cite[Chapter $8$]{slivkins2019introduction} for a detailed discussion on these last two algorithms. Finally, we also recall Thompson sampling \cite{thompson1933likelihood}. In Thompson sampling, the reward distributions are assumed to be fully described by their means, themselves sampled from a probability distribution with a finite support. Then, the history of interactions with the system is used to calculate the posterior distribution of the mean reward. This is a belief state, which is in turn used to compute an action \cite[Chapter $3$]{slivkins2019introduction}. This idea can be extended to both bandits with infinite arms (see \cite[Chapter $36$]{lattimore2020bandit} for the details) and to contextual bandits \cite{agrawal2013thompson}.\\
\newline
We close this section by making the following comments, which are transversal to our discussion on multi-armed bandit problems.\\
\newline
\noindent \textbf{Representation learning for bandits.} Inspired by humans' ability to learn and transfer experience to {\em new} tasks, a number of works have  studied representation learning \cite{bengio2013representation,representation_RL,representation_imitation} for multi-task bandits. As a prototypical multi-task DM scenario, in \cite{pasqualetti_repr} the authors consider a series of linear bandit models, where each bandit represents a different task.  The goal of the agent is to maximize the cumulative reward by interacting with these tasks. The sequential tasks are drawn from different environments, with each environment having its own representation. Given this set-up, in \cite{pasqualetti_repr} the so-called change-detection representation learning algorithm (CD-RepL) is introduced and it is shown that this algorithm outperforms certain state-of-the-art baselines. The study of how representation learning can be used to improve efficiency of bandit problems is also considered in \cite{yang2021impact}, where a set of linear bandits is studied under the assumption that an unknown linear \emph{feature extractor} exists. Both finite and infinite action settings are considered in \cite{yang2021impact} and, for each setting, an algorithm is proposed to exploit feature knowledge.\\
\newline
\noindent \textbf{Decentralized bandits.} Decentralized (or cooperative) bandits model settings where a number of agents interact with the same multi-armed bandit problem. In \cite{landgren2018social} a setting is considered where the agents are connected through an undirected graph and each agent can observe actions and rewards of its neighbors. A policy based on partitions of the communication graph is then proposed. Only one agent in each partition,  the {\em leader}, makes independent decisions based on its local information. The other agents in the partition, the {\em followers}, imitate the decisions of their leader, either directly if the leader is a neighbor, or indirectly by imitating a neighbor. In turn, the leader in each partition uses a UCB algorithm, and in \cite{landgren2018social} it is demonstrated how the group can achieve order-optimal performance. We also recall \cite{madhushani2020dynamic} where multi-agent multi-armed bandit problems in which decision-making agents can observe the choices and rewards of their neighbors are considered. Under the assumption of linear observation cost, a sampling algorithm (based on UCB) and an observation protocol are devised, which allow each agent to maximize its own expected cumulative reward.\\
\newline
\noindent \textbf{Bandits with behavior policies.} In Problem \ref{prob:multiArmed_problem} we relaxed the constraints of Problem \ref{prob:RL_problem} related to the behavior policy. While the algorithms surveyed above do not make use of a behavior policy, we highlight here a  stream of literature aimed at devising off-policy algorithms for multi-armed bandits (hence, to consider these approaches, the constraints for the behavior policy need to be added to Problem \ref{prob:multiArmed_problem}). We refer readers interested into the off-policy evaluation problem for bandits to \cite{RePEc:arx:papers:2010.12470} for a survey with a discussion on related open problems. See also   \cite{huang2021offpolicy,9741776} for a number of technical results and \cite{10.1145/3383313.3412244} for an application  to the design of a personalized treatment recommender system.
}

\section{Probabilistic design of policies through the lenses of Problem \ref{prob:general_problem}}\label{sec:control}
We consider costs that satisfy Assumption \ref{asn:cost}. That is,  we let:
\begin{equation}\label{eqn:cost_DDC_1}
  \mathbb{E}_{f}\left[c_{{1:T}}(\mathbf{X}_{0},\ldots,\bv{X}_T, \mathbf{U}_{1}, \ldots,\bv{U}_T)\right]  = \mathbb{E}_{{f}}\left[\sum_{k\in1:T}c_{k}(\mathbf{X}_{k-1}, \mathbf{U}_{k})\right],
\end{equation}
where the pf $f$ is the same as in Problem \ref{prob:general_problem} and $c_k(\cdot,\cdot)$ is the cost incurred by the agent at each $k$.  Often,  the cost in (\ref{eqn:cost_DDC_1}) is regularized and a typical choice to do so is to use some statistical divergence \cite{10.1016/j.sigpro.2012.09.003,10.2307/1403865}.  Among the possible divergences,  which offer a measure of the {discrepancy} between pairs of pfs,  a common choice is to use the so-called Kullback-Leibler  \cite{KL_51} divergence (also known as cross-entropy or relative entropy). This is formalized with the following:

\begin{Definition}\label{def:KLdiv}	
Consider two pdfs, \textcolor{black}{$\phi(\mathbf{z})$} and \textcolor{black}{$g(\mathbf{z})$,} with the former being absolutely continuous with respect to the latter. Then, the Kullback-Leibler divergence (\KL-divergence for short) of $\phi(\mathbf{z})$ with respect to $g(\mathbf{z})$, is 
	\begin{equation}
	\label{eqn:DKL_def}
	\mathcal{D}_{\KL}
	\left(\phi(\mathbf{z}) \mid \mid g(\mathbf{z}) \right):
	= \int_{\mathcal{S}(\phi)} \phi(\mathbf{z}) \; \ln\left( \frac{\phi(\mathbf{z})}{g(\mathbf{z})}\right)\,d\mathbf{z}.
	\end{equation}
\end{Definition}
Clearly,  the above definition is given for pdfs. For pmfs the same definition holds but with the integral in (\ref{eqn:DKL_def}) replaced with the sum. Given this set-up, we can now formulate the following variation on Problem \ref{prob:general_problem}:
\begin{Problem}\label{prob:DDC_problem}
Let, $\forall k\in 1:T$: (i) $\mathcal{E}_k$ and $\mathcal{I}_k$ be index sets for equality and inequality constraints; (ii) $H^{(j)}_{\bv{u},k}$,  $G^{(j)}_{\bv{u},k}$, $0\le\varepsilon_k\le 1$ be  constants; (iii) $h^{(j)}_{\bv{u},k},g^{(j)}_{\bv{u},k}:\mathcal{U}\rightarrow \R$ be  measurable mappings.  Find  $\{f^{\ast}_{\bv{u},k}(\mathbf{u}_{k}\mid \mathbf{x}_{k-1})\}_{k \in 1:T}$ such that
\begin{subequations}
    \begin{alignat}{3}
    \{f_{\bv{u},k}^{\ast}(\mathbf{u}_{k}\mid \mathbf{x}_{k-1})\}&_{k \in 1:T} \in \underset{\{f_{\bv{u},k}(\mathbf{u}_{k}\mid \mathbf{x}_{k-1})\}_{k \in 1:T}}{{\arg \min}}\  \DKL (&&f\mid \mid g) + \mathbb{E}_{{f}}\left[\sum_{k\in1:T}c_{k}(\mathbf{X}_{k-1}, \mathbf{U}_{k})\right] \label{DDCProb_cost}\\
    & s.t.  \ \E_{f_{\bv{u},k}}\left[h^{(j)}_{\bv{u},k}(\bv{U}_k)\right] = H^{(j)}_{\bv{u},k}, &&\forall k\in 1:T, \forall j \in \mathcal{E}_k;\label{DDCProb_eqConstr}\\
    & \ \ \ \ \ \E_{f_{\bv{u},k}}\left[g^{(j)}_{\bv{u},k}(\bv{U}_k)\right] \le G^{(j)}_{\bv{u},k}, &&\forall k\in 1:T, \forall j \in \mathcal{I}_k;\label{DDCProb_ineqConstr}\\
        &\ \ \ \ \ f_{\bv{u},k}(\mathbf{u}_{k}\mid \mathbf{x}_{k-1}) \in \mathcal{P}, &&\forall k\in 1:T.\label{DDCProb_pf}
    \end{alignat}
\end{subequations}
\end{Problem}
Problem \ref{prob:DDC_problem} is a regularized and relaxed version of Problem \ref{prob:general_problem} and we refer to Section \ref{sec:gen_prob_statement} for a discussion on the constraints.  The pf $f$ in the cost functional is defined in (\ref{eqn:full_joint_general}), while  the pf $g$ can be interpreted as a pf towards which the solution is biased.  In certain streams of the literature, see e.g., \cite{towardsFPDkarny,fpdKarnyGuy,mrfpdHerzallah,KARNY2020104719,KARNY2020104,QUINN2016532} and references therein, this pf is termed as the {reference} (or {ideal}) pf and expresses preferences (in terms of both performance and safety) on the desired  evolution of the closed-loop system.  In a complementary stream of literature, the pf $g$ takes the role of a passive, e.g., uncontrolled, dynamics \cite{Todorov_pnas,todorov_linearly_mdps,theodorou1,KL_as_inference}. Finally, in e.g.,  \cite{gagliardi2020probabilistic,Gagliardi_D_et_Russo_G_IFAC2020_extended_Arxiv} the reference pf is instead estimated from example data. Let $\boldsymbol{\Gamma}_{0:T}$ be the example dataset. Then, the chain rule for pfs implies that
\begin{equation}\label{eqn:full_joint_ideal}
    g:= g(\boldsymbol{\Gamma}_{0:T}) = g_0(\mathbf{x}_0)\prod_{k\in1:T}g_{\mathbf{u},k}(\mathbf{u}_{k}\mid \mathbf{x}_{k-1}) g_{\mathbf{x},k}(\mathbf{x}_{k}\mid \mathbf{x}_{k-1}, \mathbf{u}_{k}),
\end{equation}
where $g(\mathbf{x}_0)$ is a prior. In what follows we say that $g_{\mathbf{u},k}(\mathbf{u}_{k}\mid \mathbf{x}_{k-1})$ is the reference policy (e.g., extracted from the examples or embedding desired properties that the control signal should fulfill).  The pf $g_{\mathbf{x},k}(\mathbf{x}_{k}\mid\mathbf{x}_{k-1}, \mathbf{u}_{k})$ is instead termed as reference system. Such a pf can be different from  $f_{\mathbf{x},k}(\mathbf{x}_{k}\mid\mathbf{x}_{k-1}, \mathbf{u}_{k})$ in (\ref{eqn:full_joint_general}). In the context of controlling from demonstration, this means that the system used to collect the example data can be physically different from the one under control \cite{gagliardi2020probabilistic}.  Moreover,  as shown in \cite{8795901},  Problem \ref{prob:DDC_problem} is equivalent to the linear quadratic Gaussian regulator when: (i) $c_{k}(\cdot, \cdot)$ is equal to $0$ for all $k$; (ii) the only constraint is (\ref{DDCProb_pf}); (iii) all the pfs are Gaussians and the means of $g_{\mathbf{u},k}(\mathbf{u}_{k}\mid\mathbf{x}_{k-1})$ and $g_{\mathbf{x},k}(\mathbf{x}_{k}\mid\mathbf{x}_{k-1}, \mathbf{u}_{k})$ are both equal to $0$.  Interestingly,  from the information-theoretical viewpoint,  minimizing the first component in the cost functional amounts at projecting the pf $f$ onto $g$, see e.g., \cite{10.5555/1146355}.

\begin{Remark}\label{rem:classicDD1}
In Problem \ref{prob:DDC_problem} the constraint $\bv{x}_k \sim f_{\mathbf{x},k}(\bv{x}_k\mid \bv{x}_{k-1},\bv{u}_k)$ appearing in Problem \ref{prob:general_problem} is embedded in the KL-divergence component of the cost.  When $f_{\mathbf{x},k}(\bv{x}_k\mid \bv{x}_{k-1},\bv{u}_k)$ is {linear with Gaussian noise}, then this constraint is equivalent (see Section \ref{sec:pdfs}) to $\bv{X}_k = \bv{A}\bv{X}_{k-1}+\bv{B}\bv{U}_k + \bv{W}_k$,  where $ \bv{W}_k$ is sampled from a Gaussian. In this case, the constraint can be expressed in terms of the {\em behaviors} of the linear dynamics. This viewpoint is at the basis of the behavioral theory approach, pioneered by Willems starting from the $80$s \cite{4384643}.  Driven by the emergent data-driven control paradigm, behavioral theory has gained renewed interest and we refer readers to \cite{MARKOVSKY2021} for a detailed survey.  While surveying data-driven control approaches based on behavioral theory goes beyond the scope of this paper, we recall   \cite{doi:10.1080/00207170801942170,8795639,8933093,8960476,8703172,9143608} which, among others, exploit behavioral theory to tackle a broad range of data-driven control problems. 
\end{Remark}

\begin{Remark}
For nonlinear systems, another approach to represent the dynamics corresponding to $f_{\mathbf{x},k}(\bv{x}_k\mid \bv{x}_{k-1},\bv{u}_k)$ leverages the use of Koopman operators \cite{Koopman315, bruntonKoopman}. These are infinite-dimensional operators that allow to handle nonlinear systems through a globally linear representation.   See \cite{BEVANDA2021197} for a detailed survey of data-driven representations of Koopman operators for dynamical systems and for their applications to control.
\end{Remark}

\begin{Remark}\label{rem:classicDD2}
When $\bv{x}_k\sim f_{\mathbf{x},k}(\bv{x}_k\mid \bv{x}_{k-1},\bv{u}_k)$ is expressed as a difference equation, then inspiration to solve Problem \ref{prob:DDC_problem} can be gathered from MPC. Works in this direction include \cite{8039204}, where an MPC learning algorithm is introduced for iterative tasks, \cite{SALVADOR2018356} where a data-based predictive control algorithm is presented, \cite{8909368} where a MPC approach that integrates a nominal system with an additive nonlinear part of the dynamics modeled as a Gaussian process is presented. See also \cite{Recht} where coarse-ID Control methods are also surveyed.
\end{Remark}

\subsection{Finding  the optimal solution to Problem \ref{prob:DDC_problem}}

The resolution methods surveyed next, attempt to solve Problem \ref{prob:DDC_problem} by either finding a randomized policy or by finding directly a controlled transition probability, or by recasting the problem in terms of FP equations.   We refer to the first literature stream as {fully probabilistic design}, the second is termed as KL-control and the third  is the so-called FP-control.

\subsubsection{Fully Probabilistic Design}\label{FPD_section}

The stream of literature labeled as Fully Probabilistic design (FPD) tackles Problem \ref{prob:DDC_problem} when the cost only has the KL-divergence component.  In this case,  Problem \ref{prob:DDC_problem} becomes the problem of finding $\{f^{\ast}_{\bv{u},k}(\mathbf{u}_{k}\mid \mathbf{x}_{k-1})\}_{k \in 1:T}$ so that
   \begin{equation}\label{eqn:FPD_problem}
    \begin{split}
    \{f_{\bv{u},k}^{\ast}(\mathbf{u}_{k}\mid \mathbf{x}_{k-1})\}_{k \in 1:T} \in & \underset{\{f_{\bv{u},k}(\mathbf{u}_{k}\mid \mathbf{x}_{k-1})\}_{k \in 1:T}}{{\arg \min}}  \DKL (f\mid \mid g)  \\
    & \text{s.t.  constraints  (\ref{DDCProb_eqConstr}) - (\ref{DDCProb_pf})} 
    \end{split}
    \end{equation}
The above problem, termed as FPD problem in what follows, is a tracking control problem and the goal is that of designing $\{f^{\ast}_{\bv{u},k}(\mathbf{u}_{k}\mid \mathbf{x}_{k-1})\}_{k \in 1:T}$  so that the pf $f$ of the system is as similar as possible (in the KL-divergence sense) to the reference pf $g$.  To the best of our knowledge, the relaxed version of the FPD problem has been tackled for control purposes in \cite{towardsFPDkarny}. The approach builds on the Bayesian framework for system Identification \cite{bayesian_id}. See e.g.,  \cite{fpdKarnyGuy,mrfpdHerzallah,gagliardi2020probabilistic} for a set of  results that build on \cite{towardsFPDkarny}.  By assuming that state transitions can be direc{t}ly controlled (see Section \ref{sec:KL_control}) a cost has been added to the KL-divergence component. This cost can be used to formalize additional requirements on the closed-loop evolution that might not be captured by the reference pf.  For example, in e.g., \cite{crowdsourcing,crowdsourcing_regret} for a shared economy/smart cities application, the KL-divergence component models the fact that an autonomous car would like to stay on a route that accommodates the preferences of the passengers and the additional cost instead depends on road (and/or parking) conditions.   

We now proceed with surveying methods to solve the problem in (\ref{eqn:FPD_problem}). The problem admits an explicit expression for the optimal solution.   Once the optimal solution for this constrained FPD problem is presented, we then proceed with showing what happens when actuation constraints are removed.  In this case,  we find back the expression of the control policy from \cite{towardsFPDkarny}. \\

\noindent {\bf The constrained FPD.} The problem in (\ref{eqn:FPD_problem}) has been  tackled in \cite{gagliardi2020probabilistic} in the context of control synthesis from examples.  In such a paper it is shown that the problem can be broken down into convex sub-problems of the form 
\begin{equation}\label{eqn:constrained_FPD_timestep_k}	
	\begin{split}
 \underset{f_{\bv{u},k}(\bv{u}_k\mid\bv{x}_{k-1})}{\min} 
		& \mathcal{D}_{\KL}\left(f_{\bv{u},k}(\bv{u}_k\mid\bv{x}_{k-1})\mid \mid g_{\bv{u},k}(\bv{u}_k\mid\bv{x}_{k-1}) \right)+\E_{f_{\bv{u},k}}\left[\hat{\omega}\left(\mathbf{U}_k,\mathbf{X}_{k-1}\right)\right]\\
	 { \textnormal{s.t.:} }  	&  \E_{f_{\bv{u},k}}\left[h^{(j)}_{\bv{u},k}(\bv{U}_k)\right] = H^{(j)}_{\bv{u},k},  \ \ \forall j \in \mathcal{E}_k;\\ 
		&  \E_{f_{\bv{u},k}}\left[g^{(j)}_{\bv{u},k}(\bv{U}_k)\right] \le G^{(j)}_{\bv{u},k},  \ \ \forall j \in \mathcal{I}_k; \\
		&    \ f_{\bv{u},k}(\mathbf{u}_{k}\mid \mathbf{x}_{k-1}) \in \mathcal{P}.
	\end{split} 
\end{equation}
At each $k$, by solving the problem in (\ref{eqn:constrained_FPD_timestep_k}) the optimal control pf $f^{\ast}_{\bv{u},k}(\mathbf{u}_{k}\mid \mathbf{x}_{k-1})$ is obtained. In the cost functional of (\ref{eqn:constrained_FPD_timestep_k}) the term $\hat{\omega}(\cdot,\cdot)$ needs to be obtained via a backward recursion.  The results in \cite{gagliardi2020probabilistic} leverage a strong duality argument for the convex problem and require that  the constraints satisfy the following:
\begin{Assumption}\label{asm:constrraints}
There exists at least one pf that satisfies the equality constraints in Problem \ref{prob:DDC_problem} and also satisfies the inequality constraints strictly.
\end{Assumption}
The arguments in \cite{gagliardi2020probabilistic} lead to an algorithmic procedure.  The procedure takes as input $    g(\boldsymbol{\Gamma}_{0:T})$, $\{f_{\bv{x},k}(\mathbf{x}_{k}\mid \mathbf{x}_{k-1}, \mathbf{u}_{k})\}_{1:N}$ and the constraints of Problem \ref{prob:DDC_problem}.   Given this input, the algorithm outputs  
$\{f^{\ast}_{\bv{u},k}(\mathbf{u}_{k}\mid \mathbf{x}_{k-1})\}_{k\in 1:N}$ and,  at each $k$, the control input applied to the system is obtained by sampling from the pf $f^{\ast}_{\bv{u},k}(\mathbf{u}_{k}\mid \mathbf{x}_{k-1})$.  In particular, the optimal solution at time-step $k$ is given by
\begin{equation}\label{eqn:optimal_sol_constrained}
\begin{split}
&    f^\ast_{\bv{u},k}(\mathbf{u}_{k}\mid \mathbf{x}_{k-1}) = \\
&g_{\bv{u},k}(\mathbf{u}_{k}\mid \mathbf{x}_{k-1})\frac{\exp\left(-\hat{\omega}(\bv{u}_k,\bv{x}_{k-1})-\sum_{j\in\mathcal{I}_{a,k}}\lambda^{(j),\ast}_{\bv{u},k}h^{(j)}_{\bv{u},k}(\bv{u}_k)\right)}{\int g_{\bv{u},k}(\mathbf{u}_{k}\mid \mathbf{x}_{k-1})\exp\left(-\hat{\omega}(\bv{u}_k,\bv{x}_{k-1})-\sum_{j\in\mathcal{I}_{a,k}}\lambda^{(j),\ast}_{\bv{u},k}h^{(j)}_{\bv{u},k}(\bv{u}_k)\right)d\bv{u}_k},
\end{split}
\end{equation}
where $\lambda^{(j),\ast}_{\bv{u},k}$ are the Lagrange multipliers (obtained from the dual problem) and $\mathcal{I}_{a,k}$ is the set of active constraints.  While for the sake of brevity we do not report the backward recursion $\hat{\omega}(\cdot,\cdot)$,  we make the following remarks.

\begin{Remark}
Assumption \ref{asm:constrraints} becomes the classic Slater's condition when the decision variables are vectors. This assumption, in its functional form, arises in the literature on infinite-dimensional convex optimization \cite{Roc_88,Fan_K_JMAnnApp_1968_infinite, nishiyama2020convex,10.1145/2591796.2591803,Bot_05}.  The constraints are {\em moment constraints}, see e.g., \cite{1246014,1618839, 8359301}.
\end{Remark}

\begin{Remark}
The expression for the optimal solution in (\ref{eqn:optimal_sol_constrained}) defines a so-called twisted kernel \cite{BALAJI2000123}.  In the optimal solution, this twisted kernel is a Boltzmann-Gibbs distribution. Notably, in statistical physics these solutions arise as the solutions of minimization problems involving Gibbs-types free energy functionals.
\end{Remark}
{We close this paragraph by recalling  \cite{9029512,9147402,cammardella2021kullbackleiblerquadratic} where closely related infinite-dimensional finite-horizon control problems are considered in the context of distributed control and energy systems. In such papers,  control formulations are considered where state transition probabilities can be shaped directly together with the presence of an additional (quadratic) cost criterion.  See also \cite{Chertkov2018} where the minimization of a KL-divergence cost subject to moment constraints (without control variables) is considered.}\\

\noindent{\bf The unconstrained FPD.} To the best of our knowledge,  a version of the problem in (\ref{eqn:FPD_problem}) with all constraints relaxed except (\ref{DDCProb_pf}) has been originally considered in  \cite{towardsFPDkarny}.  As for the constrained case, an explicit expression for the optimal solution exists and is obtained by solving, at each $k$, a convex optimization problem.  It can be shown that the optimal solution of this unconstrained FPD problem is given by (\ref{eqn:optimal_sol_constrained}) when the functions $h_{\bv{u},k}^{(j)}$'s are all equal to $0$. Moreover,  in this case $\hat{\omega}(\cdot,\cdot)$ is obtained via the following backward recursion: 
\begin{equation}\label{eqn:optimal_FPD_backward}
\begin{split}
&\hat{\omega}(\bv{u}_k,\bv{x}_{k-1})  = \\
&\mathcal{D}_{\KL}\left(f_{\bv{x},k}(\bv{x}_k\mid\bv{x}_{k-1},\bv{u}_k)\mid \mid g_{\bv{x},k}(\bv{x}_k\mid\bv{x}_{k-1},\bv{u}_k) \right){-}\E_{f_{\bv{x},k}}\left[\ln\gamma_k(\bv{X}_k)\right] ,
\end{split}
\end{equation}
with
\begin{equation*}
\begin{split}
& \gamma_k(\bv{x}_k)  =\\
& \E_{g_{\bv{u},k+1}}\left[\exp(-\mathcal{D}_{KL}(f_{\bv{x},k+1}(\bv{x}_{k+1}\mid\bv{x}_{k},\bv{u}_{k+1})\mid \mid g_{\bv{x},k+1}(\bv{x}_{k+1}\mid\bv{x}_{k},\bv{u}_{k+1})) \right.\\
&\left. + \E_{f_{\bv{x},k+1}}\left[\ln (\gamma_{k+1}(\bv{X}_{k+1})) \right])\right], \\
&  \gamma_{T}(\bv{x}_{T}) = 1. 
\end{split}
\end{equation*}
The above solution  has been subject of algorithmic research \cite{fpdKarnyGuy} complemented with efforts towards its axiomatization \cite{KARNY2020104719,KARNY2012105}.  See also  \cite{KARNI2021,https://doi.org/10.1002/acs.743,8167323,7362235,7583697}.  For problems involving the system output (rather than the state) a solution has been presented in \cite{fpdKarnyGuy}. Finally, another line of research aimed at widening the range of conditions under which FPD can be applied considers the presence of delays. 	See e.g., \cite{fpdHerzallah,mrfpdHerzallah}, which  adapt the result to the case where, at each $k$, the dynamics of the system are conditioned on data prior to $k-1$. 

\begin{myexp*}\normalfont
The FPD framework is particularly appealing to tackle the problem of synthesizing control policies from examples.  In this case,  the reference pf $g$ is extracted from example data and captures a desired evolution that the closed loop system needs to track.  In this case, by minimizing the cost of the problem in (\ref{eqn:FPD_problem}) a randomized control policy is designed so that the pf of the closed loop system, i.e., $f$, tracks the reference pf from the examples,  i.e., $g$. That is, the policy is such that the discrepancy between $f$ and $g$ is minimized.  In this part of the running example we now make use of the pfs $g_{\bv{x}}(\bv{x}_k\mid\bv{x}_{k-1}, u_k)$ and $g_{{u}}(u_k\mid\bv{x}_{k-1})$ obtained from the reference system in the first part of the running example as reference pfs.  The pf of the system we want to control is instead $f_{\bv{x}}(\bv{x}_k\mid\bv{x}_{k-1}, u_k)$ also estimated (via Algorithm \ref{alg:hist_filter}) in the first part of the running example.  

The FPD formulation allows to tackle situations where the system under control is different from the one used to collect the example data. In fact,  in our example, the weight of the mass of the pendulum under control is different from the weight of the mass of the reference system (i.e., the pendulum under control has a mass of $1$kg, while the pendulum used to collect the example data has a mass of $0.5$kg). For concreteness,  we considered the unconstrained FPD formulation and the optimal policy  is given by (\ref{eqn:optimal_sol_constrained}) with the functions $h_{\bv{u},k}^{(j)}$'s all equal to $0$ and with $\hat{\omega}(\bv{u}_k,\bv{x}_{k-1})$ generated via (\ref{eqn:optimal_FPD_backward}).  When computing the policy we used a receding horizon strategy, with width of the horizon window $H =2$. The corresponding simulation results are reported in Figure \ref{fig:FPD_simulation}. In such a figure, it is clearly shown that the FPD policy is able to swing up the pendulum.  It can also be observed that the evolution of the controlled pendulum is similar to the one  in the examples (see Figure \ref{fig:examples_simulation}) even despite the fact that this latter pendulum is physically different from the one under control.
\begin{Remark}
The width of the receding horizon window used to obtain the FPD policy is one order of magnitude smaller than the width we used for MPC.  Still, the FPD is able to swing-up the pendulum and the intuition for this is that, if the example data are of {\em good} quality\footnote{An interesting research question is to determine  {\em minimal requirements} that make a database a good database. We refer to Section \ref{sec:conclusions} for a discussion on this important aspect.} then the control algorithm does not need to look much farther in the future in order to fulfill the control task. 
\end{Remark}
\begin{figure}[thbp]
\centering
\includegraphics[width=0.48\linewidth]{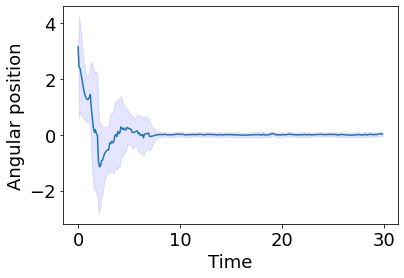}
\includegraphics[width=0.48\linewidth]{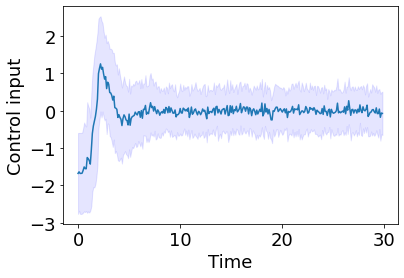}
\caption{swing-up of the pendulum via FPD.  The time evolution of $\theta_k$ is in the leftward panel, while the time evolution of $u_k$ is shown  in the rightward panel. Bold lines denote the mean profiles and the shaded areas represent the confidence intervals corresponding to the standard deviation. Figure obtained by running $50$ simulations.}
\label{fig:FPD_simulation}
\end{figure}
\end{myexp*}

\subsubsection{KL-Control}\label{sec:KL_control}
The stream of literature that goes under the label of KL-control (KLC in what follows) essentially attempts to solve the unconstrained version of  Problem \ref{prob:DDC_problem} by computing the following optimal cost-to-go function:
\begin{equation}\label{opti_costToGo}
v(\mathbf{x}_{k-1}) := \underset{\bv{u}_{k}}{\text{min\ }}l_k(\bv{x}_{k-1},\bv{u}_k) + \mathbb{E}_{f_{\bv{x},k}}[v(\bv{X}_k)].
\end{equation}
In the above expression,  $v(\mathbf{x}_{k-1}) $ is the optimal cost-to-to from state $\bv{x}_{k-1}$ and $l_k(\cdot,\cdot)$ is an immediate cost, which, as we shall see,  includes a KL-divergence component.  In (\ref{opti_costToGo}) the cost function has a component that depends on the future states and this accounts for the fact that optimal actions cannot be found, in general, via a greedy optimization of the immediate cost.  To the best of our knowledge, work on KLC can be traced back to the works by Todorov \cite{Todorov_pnas,todorov_linearly_mdps,todorov_duality_estimation,4927540,todorov_compositionality} and Kappen \cite{Kappen_path_integral_applied},  building, among others, on \cite{doi:10.1080/17442508208833228,doi:10.1137/S0363012901393894}. \\

\noindent {\bf The KLC approach to action computation.}  Actions from (\ref{opti_costToGo}) can be efficiently computed if the cost-to-go function is available.  The key idea behind KLC methods is that of finding an analytical expression for the optimal actions given the value function and then devise a transformation that linearizes (\ref{opti_costToGo}).  This approach is particularly convenient when the following assumption is made:
\begin{Assumption}\label{asn:controllability}
For each $k$:
\begin{enumerate}
    \item $f_{\bv{x},k}(\mathbf{x}_{k}\mid\mathbf{x}_{k-1},\bv{u}_k) f_{\bv{u},k}(\mathbf{u}_{k}\mid\mathbf{x}_{k-1}) = \pi(\mathbf{x}_{k}\mid \mathbf{x}_{k-1})$;
    \item $g_{\bv{x},k}(\mathbf{x}_{k}\mid\mathbf{x}_{k-1},\bv{u}_k) g_{\bv{u},k}(\mathbf{u}_{k}\mid\mathbf{x}_{k-1}) = p(\mathbf{x}_{k}\mid \mathbf{x}_{k-1})$.
\end{enumerate} 
\end{Assumption}
That is, within the KL framework, it is assumed that: (i) the agent can specify directly  the (state) transition pf rather than the action; (ii) the passive dynamics is also specified through a state transition pf.  From the control perspective, this means that the agent can specify how its state must evolve and, to do so, a low level control that guarantees the state evolution might be available. { Intuitively, one can also think of the agent's actions as probability distributions over the state space. This means that state feedback policies amount to directly choosing, for each $k$, the conditional probability law of $\bv{X}_k$ given $\bv{X}_{k-1}$. This is why there is no marginalization over $\bv{U}_k$ in Assumption \ref{asn:controllability}.}  Let $q(\cdot)$ be a state cost. Then, the immediate cost considered in the KLC framework is given by
 \begin{equation}\label{eqn:immediate_cost}
l(\bv{x}_{k-1},{\pi(\bv{x}_k\mid\bv{x}_{k-1})}) := q(\bv{x}_{k-1}) + \sD_{KL}(\pi(\mathbf{x}_{k}\mid \mathbf{x}_{k-1})\mid \mid p(\mathbf{x}_{k}\mid \mathbf{x}_{k-1})).
\end{equation}
From Definition \ref{def:KLdiv}, $\pi(\mathbf{x}_{k}\mid \mathbf{x}_{k-1})$ needs to be absolutely continuous with respect to $p(\mathbf{x}_{k}\mid \mathbf{x}_{k-1})$ and this has the interesting implication \cite{Todorov_pnas} of preventing physically impossible behaviors.  The pf $p(\mathbf{x}_{k}\mid \mathbf{x}_{k-1})$ has the role of a {\em passive} dynamics e.g., representing the evolution of the uncontrolled system.  This might define a free dynamics and deviations from this dynamics are associated to energy expenditures \cite{MDP_KL_cost}. Following \cite{Todorov_pnas}, by introducing the desirability function $z(\bv{x}_k):=\exp(-v(\bv{x}_k))$ the optimal solution to (\ref{opti_costToGo}) - (\ref{eqn:immediate_cost}) is found as 
\begin{equation}\label{eqn:optimal_sol_Todorov}
    \pi^{\ast}(\mathbf{x}_{k}\mid \mathbf{x}_{k-1}) = \frac{p(\mathbf{x}_{k}\mid \mathbf{x}_{k-1})z(\bv{x}_k)}{G_z(\bv{x}_{k-1})},
\end{equation}
where  $G_z(\bv{x}_{k-1})$ is the normalization factor.  In order to obtain the optimal solution to  (\ref{opti_costToGo}) - (\ref{eqn:immediate_cost}) the term $z(\bv{x}_k)$ needs to be computed. This can be done by plugging the analytical expression of the optimal solution into (\ref{opti_costToGo}) - (\ref{eqn:immediate_cost}). By doing so,  it can be shown that the desirability function must satisfy
\begin{equation}\label{eqn:KL_desirability}
z(\bv{x}_k) = \exp(-q(\bv{x}_k))\E_{p}\left[z(\bv{X}_{k+1})\right],
\end{equation}
which, when the state space is finite, can be recast as an eigenvector problem and in turn  this can be solved before executing the actions.

The structure of $\pi^{\ast}(\mathbf{x}_{k}\mid \mathbf{x}_{k-1})$ highlights the fact that the optimal solution twists the passive dynamics by penalizing states that are not desirable. In this context, we also recall \cite{todorov_linearly_mdps} that applies the above reasoning to controlled transition probabilities and \cite{todorov_compositionality}, where it is shown that optimal actions can be also obtained  as a sum of actions that are optimal for other problems with the same passive dynamics. It is  worth noting that \cite{todorov_duality_estimation} a duality result exists between the optimal control problem solved above and information filtering.  We also recall  \cite{MDP_KL_cost}, where the KLC framework is developed for online Markov Decision Processes with the online aspect of the problem consisting in the fact that the cost functions are generated by a dynamic environment and the agent learns the current cost only after selecting an action.  A related result is presented in \cite{Pir_21}, where, motivated by the brain's ability to reuse previous computations,  dynamic environments are considered {and the} Woodbury matrix identity is leveraged for efficient replanning of actions when the cost changes.

\begin{Remark}
The optimal solution to the KLC problem (\ref{eqn:optimal_sol_Todorov}) and to the FPD problem (\ref{eqn:optimal_sol_constrained}) have a similar structure and can be both interpreted in terms of twisted kernels. See \cite{crowdsourcing_regret} for an explicit link between the two solutions.
\end{Remark}

\begin{Remark}
From the computational viewpoint, the KLC approach relies on linearizing the Bellman equation through a nonlinear transformation. See the next part of the running example for more algorithmic details.
\end{Remark}

\begin{Remark}
As noted in \cite{gagliardi2020probabilistic}, problems of the form of (\ref{opti_costToGo}) - (\ref{eqn:immediate_cost}) and (\ref{eqn:constrained_FPD_timestep_k})  become equivalent to a maximum entropy problem when the reference pfs (or, equivalently, the passive dynamics) are uniform distributions.
\end{Remark}

\noindent {\bf Links with inference and path integrals.} As shown in \cite{KL_as_inference}, an interesting connection between (\ref{eqn:optimal_sol_Todorov}) and graphical inference exists and can be expressed through path integrals \cite{Kappen_path_integral_applied,Kappen_path_integral}. Indeed,  from the optimal pf (\ref{eqn:optimal_sol_Todorov}) we get:
$$
 \pi^{\ast}(\mathbf{x}_{1:{T}}\mid \mathbf{x}_{0}) = \frac{p(\mathbf{x}_{1:{T}}\mid \mathbf{x}_{0})\exp\left(-\sum_{k{\in 1:T}}v(\bv{x}_k)\right)}{G_z(\bv{x}_{0})},
$$
where $G_z(\bv{x}_{0})$ is the {n}ormalizing factor defined as:
\begin{equation}\label{eqn:path_integral}
G_z(\bv{x}_{0}):=\int p(\mathbf{x}_{1:{T}}\mid \mathbf{x}_{0})\exp\left(-\sum_{k{\in 1:T}}v(\bv{x}_k)\right) d\bv{x}_{1:{T}}.
\end{equation}
The shorthand notation $d\bv{x}_{1:N}$ denotes that integration in (\ref{eqn:path_integral}) is taken over the whole path of the state evolution. Also,  for notational convenience, we omit to specify that integration is taken over the domain of the pf $p(\mathbf{x}_{1:{T}}\mid \mathbf{x}_{0})$.  Note that we are considering continuous variables and an analogous definition can be obtained for discrete variables as originally done in \cite{KL_as_inference}.  It is crucial to observe that the expression in (\ref{eqn:path_integral}) is an integral over all paths (i.e., a so-called path integral \cite{4927548}) rooted from $\bv{x}_0$ and the optimal cost is given by
$-\ln G_z(\bv{x}_{0})$. With this interpretation,  the optimal solution can be obtained by computing the path integral. In \cite{KL_as_inference} it is shown that this can be done via a graphical model inference if the following two assumptions are satisfied: (i) the passive dynamics can be factorized over the components of $\bv{x}$; (ii) the matrix of interactions between components is sparse.  See \cite{KL_as_inference} for the algorithmic details.

The path integral interpretation has also been exploited in the works \cite{theodorou1,theodorou2,theodorou3,JMLR:v11:theodorou10a}, where  the method is leveraged, and further developed,  for robotics and learning systems.  See these works for a detailed theoretical and algorithmic study of path integrals in the context of control, learning and robotics. A complementary interesting idea that leverages path integrals by taking inspiration from Model Predictive Control is the so-called Model-Predictive Path Integral Control (MPPI). This algorithm, see e.g.,  \cite{7989202,doi:10.2514/1.G001921} and references therein,  finds a minimum of a KL-divergence cost by estimating the future trajectory of the system from a path integral.  We refer to e.g., \cite{9349120} which, besides presenting a  more detailed literature survey on this topic,  also introduced an extension of MPPI with robustness guarantees and demonstrates the algorithm on a real testbed.  It is also of interest to report  \cite{DBLP:conf/rss/WangSGVGLT21}, which builds on the MPPI and path integral framework to consider Tsallis divergence costs. 

\begin{myexp*}\normalfont
KLC is now used to swing-up the pendulum (with a mass of $1$kg) by finding the optimal transition pf $\pi^{\ast}(\mathbf{x}_{k}\mid \mathbf{x}_{k-1})$ solving (\ref{opti_costToGo}) - (\ref{eqn:immediate_cost}). Following the framework outlined above, the pf $p(\bv{x}_k\mid\bv{x}_{k-1})$ is the pendulum passive (i.e., uncontrolled) dynamics. In our experiments the pf was estimated via Algorithm \ref{alg:hist_filter}, from a database obtained by simulating the uncontrolled pendulum when no control was applied. Algorithm \ref{alg:hist_filter} was applied by using the same state discretization described in the first part of the running example. The state cost is instead given by $q(\bv{x}_k) := \theta_k^2 + 0.1\omega_k^2$ (recall that in the KLC framework the state cost encourages the agent to depart from its rest position captured via the passive dynamics).

A key algorithmic feature of KLC, which makes it particularly appealing to tackle sequential decision-making problems involving transition pfs, comes from observing that (\ref{eqn:KL_desirability}) is linear in $z(\cdot)$. Hence, once the state space is discretized, the states can be enumerated (say, from $1$ to $s$) and one can then represent $z(\cdot)$ and $q(\cdot)$ as vectors, say $\bv{z}$ and $\bv{q}$. That is, the equality in (\ref{eqn:KL_desirability}) becomes
\begin{equation}\label{eqn:z_iteration}
\bv{z} = \text{diag}(\exp(-\bv{q}))\bv{P} \bv{z},
\end{equation}
where $\text{diag}(\exp(-\bv{q}))$ is the diagonal matrix having on its main diagonal the elements $\exp(-{q}(\bv{x}_1)), \ldots, \exp(-{q}(\bv{x}_s))$ and $P$ is the $s \times s$ matrix having as element $(i,j)$ the probability of transitioning from state $\bv{x}_i$ to $\bv{x}_j$. Given this set-up, computing the desirability vector $\bv{z}$ amounts at solving an eigenvector problem. We used the power iteration method to solve (\ref{eqn:z_iteration}) and hence find $\bv{z}$. Once this was obtained, then the optimal transition pf was computed via (\ref{eqn:optimal_sol_Todorov}). 

The effectiveness of KLC, implemented via the process outlined above, is shown in Figure \ref{fig:KL_simulation}. Such a figure clearly shows that the optimal transition pf $\pi^{\ast}(\mathbf{x}_{k}\mid \mathbf{x}_{k-1})$ effectively swings-up the pendulum, stabilizing the unstable upward equilibrium. From the figure, we note the following: (i) the standard deviation is smaller than the one observed in the FPD experiments. Indeed, we observed that the desirability function twisted the transition pfs to make them concentrated in a few states. The reduced standard deviation when compared to the numerical results via FPD can be  explained by the fact that FPD returns a policy that is randomized (KLC instead does not return policies but transition pfs); (ii) the behavior of the pendulum (first performing a partial swing up in the positive angles, before reaching the upward position through negative angles) can be explained in terms of the interplay between the passive dynamics and the state cost, which express{es} two different goals for the agent.

\begin{figure}[thbp]
\centering
\includegraphics[width=0.9\linewidth]{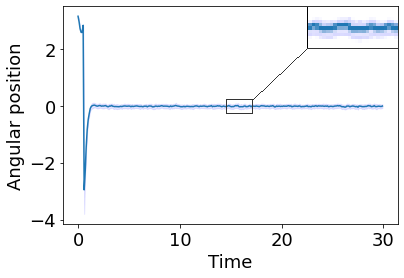}
\caption{swing-up of the pendulum via KLC. Bold line is the mean profile, the shaded area represents the confidence interval corresponding to the standard deviation. Figure obtained from $50$ simulations (the zoom magnifies the confidence interval).}
\label{fig:KL_simulation}
\end{figure}
\end{myexp*}

\subsubsection{The Fokker-Planck control framework}\label{sec:FP_control}

For completeness, we also report  an alternative approach to solve Problem \ref{prob:DDC_problem} that goes under the label of FP-control (FPC in what follows). The key idea of FPC is that of tackling Problem \ref{prob:DDC_problem} by recasting the process that is generating the data as a FP equation (see  Section \ref{sec:pdfs}).  The  framework surveyed in this section has been considerably developed in a number of works, including \cite{FP1,10.1016/j.cam.2012.06.019,quantumControl,Ann_Bor_Nob_Tem_14} and references therein. As noted in \cite{FP_survey}, in order to formulate the FPC problem the following key ingredients are needed:
\begin{enumerate}
\item the definition of a control function,  say $\bv{u}_t:=\bv{u}(\bv{x},t)$, that drives the stochastic process;
\item a FP reformulation of the controlled stochastic process;
\item the cost functional. 
\end{enumerate}

The framework allows to consider problems with actuation constraints and these are formally expressed by imposing that the function $\bv{u}_t$ belongs to a closed and convex set of admissible inputs, say $\bv{u}_t\in\mathcal{U}\subseteq\mathcal{U}_H$,  where $\mathcal{U}_H$ is a Hilbert space.   For the developments of the theory,  it is  assumed that the pf of the controlled stochastic process also belongs to a Hilbert space.  Given this set-up,  the problem of determining a control function $\bv{u}_t$ so that, starting from some initial distribution, the closed-loop system evolves, at time $t=T$, towards a desired probability density $\rho_d(\bv{x},t)$ can be formalized as follows
\begin{subequations}\label{eqn:FP_problem}	
	\begin{alignat}{3} 
 \underset{\bv{u}_t\in\mathcal{U}}{\min} & \frac{1}{2}\Vert \rho(\cdot ,T)-\rho_d(\cdot, T)\Vert_{L^2(\mathcal{X})}^2 + \frac{\nu}{2}\Vert u \Vert_{L_2(\mathcal{X} \times (0,T))}^2 \\
	 { \textnormal{s.t.:} }  	& \partial_t \rho(\bv{x},t) + \sum_{i{\in 1:n_x}} \partial_{x_i}\left(b_i(\bv{x},t,\bv{u})\rho(\bv{x},t)\right) - \sum_{i,j{\in 1:n_x}}\partial^2_{x_ix_j}\left(a_{ij}(\bv{x},t)\rho(\bv{x},t)\right) = 0,\\
	 & \rho(\bv{x},0) = \rho_0(\bv{x}).
	\end{alignat} 
\end{subequations}
In the above expression,  $\nu \ge 0$ and the constraints are the FP equation obtained from the controlled SDE
$$
d\bv{X}_t = b(\bv{X}_t,t,\bv{u}_t)dt + \sigma(\bv{X}_t,t)d\bv{W}_t.
$$
The coefficients in the constraint of the  problem in (\ref{eqn:FP_problem}) are defined as in Section \ref{sec:pdfs}. Also,  the subscript $L^2(\cdot)$ denotes the $L^2$ norm over the sets included in its parentheses.  Hence, the cost functional captures the fact that the distance between the pf of the closed-loop system and the reference pf is minimized at time $T$. The optimization problem, which can be solved via first-order necessary optimality conditions \cite{FP_survey},  has intrinsically a continuous-time formulation. However,  the problem can be discretized and this approach leads to the recursive schemes  introduced in e.g., \cite{FP1,FP2}.  The schemes rely on: 
\begin{enumerate}
\item discretizing the time domain into intervals, say $(t_k,t_{k+1})$ with $t_0 = 0$, $t_N = T$ and $t_k< t_{k+1}$; 
\item solving the optimization problem in (\ref{eqn:FP_problem}) in each interval $(t_k,t_{k+1})$ setting the initial pf $f(\bv{x},t_k) := \rho_k(\bv{x})$; 
\item keeping, in each $(t_k,t_{k+1})$, the control input $\bv{u}_k$ set to the optimal solution of the optimization problem in that time interval. 
\end{enumerate}
Within this scheme, inspired by MPC, in each time interval the optimization problem can be solved either numerically or analytically (for example, \cite{FP1} demonstrates the use of a finite elements method) to find the optimal input $\bv{u}_k$ in the time interval $[t_k, t_{k+1}]$.

Here we report a number of  works that build on the FP approach. We recall \cite{FP_hamiltonian}, which replaces the cost in (\ref{eqn:FP_problem}) with a joint cost consisting of a first term containing the expectation of a given function of the state at the end of the time interval and of a second term that is the integral of an expected cost. In both terms, expectations are taken over  the pf $\rho(\bv{x},t)$ and the problem is solved using a sampling-based Hamiltonian estimation. Moreover, in \cite{pmpc}, the system's evolution is assumed to be fully deterministic and this leads to a  FP equation that leverages Monte-Carlo samples. In \cite{simplifiedCost} a polynomial control strategy is proposed to minimize a cost that depends on the derivatives of the state's probability potentials, while \cite{forbes2004decomposition} uses a Gram-Charlier parametrization on the reference pf and calculates the control law by injecting this parametrization into a stationarity condition on the state's pf.  Finally,  \cite{stationaryFP} builds on this principle by calculating the FP equation for this stationary pf.

\begin{Remark}
As surveyed in e.g.,  \cite{FP_survey},  FPC  can be extended to consider any process that can be recast as FP equation. These processes go beyond the SDEs considered in Section \ref{sec:pdfs} and include diffusions with jumps as well as piecewise-deterministic processes.  Also,  by leveraging mean field arguments, the FP approach has proved to be an effective tool to control large-scale systems and to study what happens to these systems as the number of agents increases.
\end{Remark}

\section{Concluding discussion}\label{sec:conclusions}

We discussed the relevance of sequential DM problems that involve optimizing over probability functions for learning and control. The survey was organized around a framework that consists of the formulation of an infinite-dimensional sequential DM problem and of a set of methods to search through probability functions.   Problem \ref{prob:general_problem} served as an overarching formulation through which we revisited a wide range of popular learning and control algorithms.  This was done by proposing suitable variations on the problem and by subsequently exploring different resolution methods to find the optimal solution for these variations. We used a running example to complement our discussion. From the survey, a number of key challenges naturally arise. These are discussed next.  

From the methodological viewpoint, a first key challenge is to extract useful knowledge (in terms of e.g., pfs) from data when running experiments comes at a cost.  A common feature of learning and control approaches that rely solely on the available data  is that these need to be sufficiently informative.  Ideas from data-informativity \cite{8960476,8933093,COLIN2020109000} and optimal experiment design \cite{9406124,10.2307/2238570} might be leveraged to define a metric quantifying the {\em value of information} (similar in spirit to e.g., \cite{9541011}, see also \cite{doi:10.1287/educ.1080.0039}) gained at the expenses of new experiments.  Another challenge is the design of decision-making mechanisms that are able to tackle new situations that have never been seen before by the decision-maker.  The ability of answering these {\em what if} questions is typical of agents that can reason counterfactually. {While decision-making techniques relying on running {\em what if} type simulations are available (e.g., in simulation-based control with rollouts and MPC) a principled synthesis of control techniques with modern ideas of counterfactual causality science \cite{10.5555/1642718,Scholkopfetal21} appears to be a pressing open challenge.} Another challenge arises from the fact that, nowadays, objects are becoming smaller, smarter and with the ability of being interconnected: in one word, objects are now {\em data-sources}. In this context, a key challenge is to design agents  able to make decisions by re-using  this distributed information, without having to necessarily gather new data.  A way to tackle this challenge might be the design of open and transparent crowdsourcing mechanisms for autonomous decision-makers  \cite{crowdsourcing,crowdsourcing_regret}. These mechanisms, besides giving performance guarantees should also be able to handle cooperative/competitive behaviors among peers \citep{Anastassacos_Hailes_Musolesi_2020}. Further,  studies in neuroscience (see \cite{1000brains} for an introductory review of the theory, which is based on the pioneering work \cite{mountcastle1978organizing}) hint that similar crowdsourcing mechanisms might be implemented by the brain's {neocortex} to orchestrate how the output of certain cortical circuits are used to build models and actions. In turn, it is believed that this mechanism might be at the basis of our ability to re-use acquired knowledge in order to synthesize new action policies tackling increasingly complex tasks.

From the application viewpoint,  we see as a pressing open challenge that of establishing a widely accepted set of metrics across control and learning.  Much effort has been devoted to create a {suite} of in-silico environments and datasets \cite{brockman2016openai,6386109,D4RL,antonova2021dynamic} to test learning and control algorithms. However, these algorithms are often benchmarked only via their reward diagrams. For control applications in physical environments, these diagrams should be complemented not only with a  set of metrics that quantifies typical control performance as a function of the data available to the decision-maker (a first effort in this direction can be found in \cite{delellis2020tutoring,delellis2021controltutored}) but also with a metric that quantifies the energy consumption needed to compute the policy. In this context, we see the approximation of the policies via probabilistic graphical models and neuromorphic computing as a promising approach \citep{9612647}. The computational aspect (both for learning and control) is a challenge in its own and different methods need to be rigorously benchmarked along this dimension.  Approximation results exist that allow to reduce the computational burden and we highlight a perhaps less explored direction to integrate data-driven and model-based {\em technologies} so that they tutor each other \cite{delellis2020tutoring,9610789,10.1007/978-3-030-55180-3_6}. Finally, we believe that
the ultimate challenge will be to deploy the algorithms underpinned by the methods presented here in applications where reliable models are intrinsically probabilistic and/or hard/expensive to find.  We see quantum computing  \cite{quantumControl,quantum_ML},   biochemical systems \cite{cellMigration,fang2021stochastic}, learning/control applications with humans in the loop \cite{zhan2021human,lien2009interactive,8357977} {and the design of autonomous agents reliably executing tasks in unknown, non-stationary and stochastic environments,}  as potential test-beds that are particularly appealing for the methods surveyed here.  We hope that the framework we presented will contribute to map which method is best-suit{ed} for each of the application areas.



\end{document}